\documentclass[12pt]{amsart}
\usepackage{amsmath,amsfonts,amssymb,latexsym}
\usepackage{hyperref}
\usepackage{graphicx}
\usepackage[a4paper, left=3cm, right=3cm, top=3cm, bottom=3cm]{geometry}

\newcommand{\NN}{\mathbb{N}}
\newcommand{\CC}{\mathbb{C}}
\newcommand{\RR}{\mathbb{R}}

\newcommand{\EE}{\mathbb{E}}
\newcommand{\Oo}{\mathcal{O}}
\newcommand{\Bo}{\mathcal{B}}
\newcommand{\RE}{ {\rm Re \,} }

\newtheorem{Thm}{Theorem}
\newtheorem{Lem}{Lemma}
\newtheorem{Prop}{Proposition}
\newtheorem{Cor}{Corollary}
\theoremstyle{definition}
\newtheorem{Rem}{Remark}
\newtheorem{Def}{Definition}
\newtheorem{Ex}{Example}

\begin{document}
\title{On sequences preserving summability}
\author{S{\l}awomir Michalik}
\address{Faculty of Mathematics and Natural Sciences,
College of Science,
Cardinal Stefan Wyszy{\'n}ski University,
W\'oycickiego 1/3,
01-938 Warszawa, Poland,
ORCiD: 0000-0003-4045-9548}
\email{s.michalik@uksw.edu.pl}
\author{Maria Suwi{\'n}ska}
\address{Faculty of Mathematics and Natural Sciences,
College of Science,
Cardinal Stefan Wyszy{\'n}ski University,
W\'oycickiego 1/3,
01-938 Warszawa, Poland,
ORCiD: 0000-0002-4493-6519}
\email{m.suwinska@op.pl}
\author{Bożena Tkacz}
\address{Faculty of Mathematics and Natural Sciences,
College of Science,
Cardinal Stefan Wyszy{\'n}ski University,
W\'oycickiego 1/3,
01-938 Warszawa, Poland,
ORCiD: 0000-0002-0099-0871}
\email{bpodhajecka@o2.pl}

\begin{abstract}
 In this paper, we study sequences of positive numbers preserving summability. In particular, the open set property for such a family of sequences is shown. Several classes of sequences preserving summability, including polynomials, sums of powers, and solutions of some difference equations, are introduced.
\end{abstract}
\maketitle
\section{Introduction}
The concept of sequences preserving summability was introduced fairly recently in~\cite{I-M} as a way to study summable power series solutions of some linear $q$-difference-differential equations. However, it also happened to be very fruitful in a more general framework of summable formal power series solutions of moment differential equations. In particular, it allows to reduce the study of summable solutions of an entire family of moment differential equations to the study of such solutions for a single given moment differential equation. 

A parallel idea of sequences preserving $q$-Gevrey asymptotic expansions was also studied in~\cite{L-M}. It was shown in \cite{L-M} that such family of sequences is contained within the family of sequences preserving summability. Moreover, it was posed (see \cite[Introduction]{L-M}) as an open question whether this inclusion is strict or not. 

The present work is devoted to further study of the theory of sequences of positive real numbers that preserve summability. Specifically, we prove that the family of sequences preserving summability has the following open set property: 

\begin{quote}
    if $\tilde{m}$ is a sequence preserving summability and $m$ is a sequence close to $\tilde{m}$ in some sense, then also $m$ is a sequence preserving summability.
\end{quote}

Using both this property and a similar trait of sequences preserving $q$-Gevrey asymptotic, we construct a sequence preserving summability, which however does not preserve $q$-Gevrey asymptotic expansion. From this we conclude that the inclusion considered in~\cite{L-M} is strict, which in turn solves the open
problem posed in~\cite{L-M}.

Moreover, we describe new classes of sequences of positive real numbers that preserve summability. In particular, we show that positive polynomials and sums of exponents are in fact sequences preserving summability. We also prove that more general sequences $(m(n))_{n\geq 0}$ of positive numbers, which are solutions of homogeneous difference equations with constant coefficients, preserve summability as well. To prove these results we use the special functional equations, which define the analytic continuation of a given function together with the estimation of its growth at infinity. This way, we obtain the analytic continuation of the sum of series given in the form $\sum\frac{t^n}{m(n)}$ to the function that is analytic on the set $\CC\setminus\RR_+$ with a less than exponential growth at infinity.
In the authors' opinion, this result is intrinsically interesting.

The paper is organized as follows. In Section~\ref{sec:2} we recall necessary definitions and previously known results linked with sequences preserving summability and sequences preserving $q$-Gevrey asymptotic expansions. In Section~\ref{sec:3} we provide the proof of the fact that the family of sequences preserving summability and the family of sequences preserving $q$-Gevrey asymptotic expansions have the open set property (Theorems~\ref{th:3} and~\ref{th:5}).
In this section we also solve, in Theorem~\ref{th:7}, the open problem posed in~\cite{L-M}. Next, in Section~\ref{sec:4} we find several classes of sequences of positive numbers preserving summability, such as: polynomials (Theorem~\ref{th:th_8}), sums of powers (Theorem~\ref{th:9}) and solutions of homogeneous difference equations with constant coefficients (Theorem~\ref{th:10}). Finally, in Section~\ref{sec:5}, we discuss the application of sequences preserving summability to the solutions of moment differential equations. We also pose the conjecture that the set of sequences preserving summability is not closed under addition. 

\section{Preliminaries}\label{sec:2}
\subsection{Notation}
Let $\EE$ denote a complex Banach space with a norm $\|\cdot\|_{\EE}$. The vector space of formal power series with coefficients in $\EE$ is denoted by $\EE[[t]]$.
For any given set $G\subseteq\CC$, by $\Oo(G,\EE)$ we denote the set of all $\EE$-valued holomorphic functions defined on some open set $U$ containing $G$. We also write $\Oo(G)$ if $\EE=\CC$ for simplicity.

An unbounded sector $S$ in a direction $d\in\RR$ with an opening $\alpha>0$ is defined by
\begin{displaymath}
S=S_d(\alpha):=\{z\in\CC\colon\ z=r\*e^{i\phi},\ r>0,\ \phi\in(d-\alpha/2,d+\alpha/2)\}.
\end{displaymath}
If the opening $\alpha$ is not essential,
the sector $S_d(\alpha)$ is briefly denoted by $S_d$.
 
 An open disc in $\CC$ of radius $r>0$ with a center at the origin is denoted by $D_r$. 
 In case when the radius $r$ is not essential, the set $D_r$ is briefly denoted by $D$. We also briefly denote a disc-sector $S_d(\alpha)\cup D_r$ by $\hat{S}_d(\alpha,r)$ (resp. $S_d\cup D$ by $\hat{S}_d$).
\subsection{Sequences and formal series}
Throughout the paper, if we do not specify otherwise, it is assumed that $m=(m(n))_{n\geq 0}$ (resp. $\tilde{m}=(\tilde{m}(n))_{n\geq 0}$) is a sequence of real positive numbers with $m(0)=1$ (resp. $\tilde{m}(0)=1$). 

First, we introduce some essential definitions related to sequences and formal power series.
\begin{Def}
\label{df:order}
Let $s\in\RR$ and let $m=(m(n))_{n\geq 0}$ be a sequence. If there exist constants $0<a<A<\infty$ such that $a^n(n!)^s\leq m(n) \leq A^n (n!)^s$ for every $n\in\NN_0$ then $m$ is called a \emph{sequence of order $s$}.
\end{Def}

\begin{Ex}
Let $\Gamma(\cdot)$ be the gamma function. For fixed $k>0$,
the sequence of positive numbers $\Gamma_{1/k}:=(\Gamma(1+n/k))_{n\ge0}$ is a sequence of order $1/k$.
\end{Ex}

\begin{Def}[see {\cite[Section 5.2]{B2}}]
\label{df:borel}
 For a fixed sequence $m=(m(n))_{n\geq 0}$ of positive numbers, a linear operator $\Bo_{m}\colon \EE[[t]]\to\EE[[t]]$ defined by
 \begin{equation*}
  (\Bo_{m}\hat{u})(t):=
  \sum_{n=0}^{\infty}\frac{a_n}{m(n)}t^n\quad\text{for}\quad\hat{u}(t)=\sum_{n=0}^{\infty}a_nt^n\in\EE[[t]]
 \end{equation*}
 is called an \emph{$m$-Borel operator}.
\end{Def}

\begin{Def}
    \label{df:gevrey}
    Let $s\in\RR$. A series
    $\hat{u}(t)=\sum_{n=0}^{\infty}a_nt^n\in\EE[[t]]$ is called a \emph{formal power series of Gevrey order $s$} if
    there exist constants $B,C<\infty$ such that
    \begin{equation}
    \label{eq:est_un}
     \|a_n\|_{\EE}\leq BC^n(n!)^s\quad\textrm{for every}\quad n\in\NN_0.
    \end{equation}
    The space of formal power series of Gevrey order $s$ is denoted by $\EE[[t]]_s$.
\end{Def}

\begin{Rem}
    \label{re:gevrey}
    Let $s\in\RR$ and $m$ be a sequence of order $s$. By the definitions listed above, $\hat{u}\in\EE[[t]]_s$ if and only if there exists a disc $D$ with a center at the origin such that ${\Bo}_{m}\hat{u}\in\Oo(D,\EE)$.
\end{Rem}

\begin{Def}
 A series $\hat{u}(t)\in\EE[[t]]$ is called a \emph{formal power series of Gevrey order $-\infty$} if $\hat{u}(t)$ is a series of Gevrey order $-1/k$ for every $k>0$.
 
 The space of formal power series of Gevrey 
    order $-\infty$ is denoted by 
 \begin{equation*}
  \EE[[t]]_{-\infty}=\bigcap_{k>0}\EE[[t]]_{-1/k}.
 \end{equation*}
\end{Def}

\subsection{Growth of functions}
Functions of exponential and less than exponential growth will play a crucial role in our study.
 
\begin{Def}
\label{df:growth}
A function $f\in\Oo(\hat{S}_d(\alpha,r),\EE)$
is of \emph{exponential growth of order at most $k>0$ as $t\to\infty$
in $\hat{S}_d(\alpha,r)$} if for any
$\tilde{\alpha}\in(0,\alpha)$ and $\tilde{r}\in(0,r)$ there exist positive constants
$A,B$ such that
\begin{equation*}
\|f(t)\|_{\EE}<Ae^{B|t|^k}\quad \textrm{for every}\ t\in \hat{S}_d(\tilde{\alpha},\tilde{r}).
\end{equation*}
The space of such functions is denoted by $\Oo^k(\hat{S}_d(\alpha,r),\EE)$.
\end{Def}

\begin{Def}
We say that a function $f\in\Oo(\hat{S}_d,\EE)$
is of \emph{less than exponential growth as $t\to\infty$
in $\hat{S}_d$} if $f\in\Oo^k(\hat{S}_d,\EE)$ for every $k>0$.

We will denote the space of such functions by $\Oo^{>0}(\hat{S}_d,\EE)$.
\end{Def}

\begin{Rem}
 Observe that for $\alpha>2\pi$ we receive $\hat{S}_d(\alpha,r)=\CC$. Hence, the definitions listed above include also entire functions of exponential growth and of less than exponential growth. 
\end{Rem}

By the Cauchy estimates for derivatives we receive the following 
\begin{Prop}
\label{pr:1}
Let us assume that $f$ is an entire function and $\hat{f}$ denotes its Taylor expansion, i.e.,
$$
\hat{f}(t)=\sum_{n=0}^{\infty}f_nt^n=\sum_{n=0}^{\infty}\frac{f^{(n)}(0)}{n!}t^n.
$$
Then for every $k>0$
\begin{equation*}
 f\in\Oo^{k}(\CC,\EE)\quad\text{if and only if}\quad\hat{f}\in\EE[[t]]_{-1/k}.
\end{equation*}
Moreover,
\begin{equation*}
 f\in\Oo^{>0}(\CC,\EE)\quad\text{if and only if}\quad\hat{f}\in\EE[[t]]_{-\infty}.
\end{equation*}
\end{Prop}

\subsection{Sequences preserving summability}
In this section we introduce $k$-summable power series and the concept of sequences preserving summability. For more details about summability we refer the reader to \cite{B2,LR}.
\begin{Def}
\label{df:summable}
Let $k>0$ and $d\in\RR$. Then $\hat{u}\in\EE[[t]]$ is called
\emph{$k$-summable in a direction $d$} if there exists a disc-sector $\hat{S}_d$ in a direction $d$ such that
$\Bo_{\Gamma_{1/k}}\hat{u}\in\Oo^k(\hat{S}_d,\EE)$.

The space of $k$-summable formal power series in a direction $d$ is denoted by $\EE\{t\}_{k,d}$.
\end{Def}
\begin{Def}[{\cite[Definition 11]{I-M}}]
\label{df:preserve}
  We say that a sequence $m=(m(n))_{n\geq 0}$ \emph{preserves summability} if for any $k>0$, $d\in\RR$ and any $\hat{u}\in\EE[[t]]$ the following equivalence holds:
  \begin{equation*}
   \hat{u}\in\EE\{t\}_{k,d}\quad\text{if and only if}\quad
   \Bo_{m}\hat{u}\in\EE\{t\}_{k,d}.
  \end{equation*}
\end{Def}

\begin{Rem}
\label{re:1}
By \cite[Remark 9]{I-M} every sequence preserving summability is a sequence of order zero.
\end{Rem}

\begin{Ex}[{\cite[Example 3]{I-M}}]
 Not every sequence of order zero preserves summability. For example the sequence
 \begin{equation*}
  m(n):=
    \left\{
    \begin{array}{ll}
     1,& \textrm{if }n\ \text{is even}\\
     2^{-1}& \textrm{if } n\ \text{is odd}
    \end{array}
    \right.,
 \end{equation*}
 does not preserve summability.
\end{Ex}

\begin{Rem}
\label{re:groups}
 The set of sequences preserving summability forms a group with a group operation given by multiplication. If $m_1=(m_1(n))_{n\geq0}$ and $m_2=(m_2(n))_{n\geq 0}$ preserve summability then also their product $m=m_1\cdot m_2$ (that is $m=(m(n))_{n\geq 0}$, where $m(n)=m_1(n)\cdot m_2(n)$ for any $n\in\NN_0$) preserves summability. Observe also that the identity element $\mathbf{1}=(1)_{n\geq 0}$ and the inverse element $m^{-1}=(m(n)^{-1})_{n\geq 0}$ to $m=(m(n))_{n\geq 0}$ preserve summability.
\end{Rem}

Directly from Definitions~\ref{df:borel} and~\ref{df:preserve} we get
\begin{Prop}
Let $m=(m(n))_{n\geq 0}$ be a sequence preserving summability. Then the Borel transform $\Bo_m$ is a linear automorphism with the inversion given by $(\Bo_m)^{-1}=\Bo_{m^{-1}}$ on the following spaces of power series:
\begin{enumerate}
    \item[a)] the space $\EE[[t]]$ of formal powers series,
    \item[b)] the space $\EE[[t]]_s$ of formal power series of Gevrey order $s$ for every $s\in[-\infty,+\infty)$,
    \item[c)] the space $\EE\{t\}_{k,d}$ of $k$-summable formal power series in a direction $d$ for every $k>0$ and $d\in\RR$.
\end{enumerate}
\end{Prop}

We have the following characterization of sequences preserving summability. 
 \begin{Thm}[{\cite[Theorem 1]{I-M}}]
 \label{th:1}
  A sequence $m=(m(n))_{n\geq 0}$ preserves summability if and only if for every $k>0$ and for every $\theta\neq 0\mod 2\pi$ there exists a disc-sector $\hat{S}_{\theta}$ such that
  \begin{equation}
  \label{eq:eq_1}
  \Bo_{m}\big(\sum_{n=0}^{\infty}t^n\big)\in\Oo^k(\hat{S}_{\theta})\quad \textrm{and}\quad \Bo_{m^{-1}}\big(\sum_{n=0}^{\infty}t^n\big)\in\Oo^k(\hat{S}_{\theta}).
  \end{equation}
 \end{Thm}
 \begin{Rem}
 \label{re:4}
  Since the condition $f\in\Oo^{>0}(\CC\setminus\RR_+)$ means specifically that for every $k>0$ and for every $\theta\neq 0\mod 2\pi$ there exists a disc-sector $\hat{S}_{\theta}$ such that $f\in\Oo^k(\hat{S}_{\theta})$, by Theorem \ref{th:3} a sequence $m=(m(n))_{n\geq 0}$ preserves summability if and only if
  \begin{equation}
  \label{eq:eq_1.5}
  \Bo_{m}\big(\sum_{n=0}^{\infty}t^n\big)\in\Oo^{>0}(\CC\setminus\RR_+)\quad \textrm{and}\quad \Bo_{m^{-1}}\big(\sum_{n=0}^{\infty}t^n\big)\in\Oo^{>0}(\CC\setminus\RR_+).
  \end{equation}
 \end{Rem}

Let us now present several known examples of sequences preserving summability.

\begin{Ex}[{\cite[Example 2]{I-M}}]\label{ex:4}
If $a>0$ and $\mathbf{a}:=(a^n)_{n\geq 0}$ then the sequence $\mathbf{a}$ preserves summability. In particular, the sequence $\mathbf{1}=(1)_{n\geq 0}$ preserves summability in a trivial way.
\end{Ex}

\begin{Ex}[{\cite[Example 2]{I-M}}]
By Balser's theory of general summability \cite[Section~6.5 and Theorem~38]{B2} for any moment function $\mathfrak{m}(u)$ of order zero (see \cite[Definition~4]{Mic8} for a definition of a moment function of order zero), the sequence $(\mathfrak{m}(n))_{n\geq 0}$ preserves summability.
In particular, if
\begin{equation*}
\mathfrak{m}(n):=\frac{\Gamma(1+a_1 n)\cdots\Gamma(1+a_k n)}{\Gamma(1+b_1 n)\cdots \Gamma(1+b_l n)}\quad \text{for}\quad n\in\NN_0,
\end{equation*}
where $a_1,\dots,a_k$ and $b_1,\dots b_l$ are positive numbers satisfying
\begin{equation*}
 a_1+\cdots+a_k=b_1+\cdots+b_l,
\end{equation*}
then the sequence $(\mathfrak{m}(n))_{n\geq 0}$ preserves summability.
\end{Ex}

\begin{Ex}[{\cite[Theorem 2]{I-M}}]\label{ex:5}
For every $n\in\NN_{0}$ we define a $q$-analog of $n$ by 
\begin{equation*}
 [n]_q:=1+q+\dots+q^{n-1}=\frac{1-q^n}{1-q}.
\end{equation*}
We also introduce a $q$-analog of the factorial $n!$
\begin{equation*}
 [n]_q!:=
    \left\{
    \begin{array}{ll}
     1&\textrm{for}\ n=0\\  
     {[1]_q\cdots [n]_q}&\text{for}\ n\geq 1
     \end{array}
    \right..
\end{equation*}
If $q\in[0,1)$ then the sequence $([n]_q!)_{n\geq 0}$ preserves summability.
\end{Ex}

\subsection{Sequences preserving $q$-Gevrey asymptotic expansions}
In this section we recall the concept of sequences preserving $q$-Gevrey asymptotic expansions introduced in \cite{L-M} in a similar way to sequences preserving summability.
\begin{Def}
Let $s\in\RR$. The set of $(q;s)$-Gevrey formal power series, denoted by $\mathbb{E}[[t]]_{q;s}$, consists of all formal power series $\hat{u}(t)=\sum_{n=0}^{\infty}a_nt^n\in\mathbb{E}[[t]]$ such that there exist constants $A, B>0$, and $\alpha\in\RR$, for which:
$$\left\|a_n\right\|_{\mathbb{E}}\le A B^{n}(n!)^{\alpha}q^{s\frac{n(n-1)}{2}},\qquad n\in\NN_0.$$
We also define the set of $(q;-\infty)$-Gevrey formal power series, denoted by $\mathbb{E}[[t]]_{q;-\infty}$, as the set of all formal power series, which are $(q;s)$-Gevrey series for every $s<0$. In other words
\begin{equation*}
 \mathbb{E}[[t]]_{q;-\infty}:=\bigcap_{s<0}\mathbb{E}[[t]]_{q;s}.
\end{equation*}
\end{Def}

\begin{Def}
Let $\hat{u}(t)=\sum_{n\ge0}a_nt^n\in\mathbb{E}[[t]]$. The \emph{formal $q$-Borel transformation of order $s>0$ of $\hat{u}$} is defined by
$$\mathcal{B}_{q;s}(\hat{u})(t)=\sum_{n\ge0}\frac{a_n}{q^{s\frac{n(n-1)}{2}}}t^n.$$
\end{Def}

\begin{Def}
Let $d\in\RR$ and $s>0$. Given $f\in\mathcal{O}(\hat{S}_d,\mathbb{E})$, we say that $f$ is \emph{of $q$-exponential growth of order $1/s$} if there exist $C,h>0$ and $\alpha\in\RR$ such that
\begin{equation}
\left\|f(t)\right\|_{\mathbb{E}}\le C\exp\left(\frac{\log^2(|t|+h)}{2s\log(q)}\right)(|t|+h)^{\alpha}\quad\text{for every}\ t\in \hat{S}_d.
\end{equation}
We will denote the space of such functions by $\mathcal{O}^{q;1/s}(\hat{S}_d,\mathbb{E})$.
\end{Def}

\begin{Def}
We say that a function $f\in\Oo(\hat{S}_d,\EE)$
is of \emph{less than $q$-exponential growth as $t\to\infty$
in $\hat{S}_d$} if $f\in\Oo^{q;k}(\hat{S}_d,\EE)$ for every $k>0$.

The space of such functions is denoted by $\Oo^{q;>0}(\hat{S}_d,\EE)$, i.e.,
\begin{equation*}
 \Oo^{q;>0}(\hat{S}_d,\EE)=\bigcap_{k>0}\Oo^{q;k}(\hat{S}_d,\EE).
\end{equation*}
\end{Def}
\begin{Rem}
 For every $k_2>k_1>0$ and $l_2>l_1>0$ we have the proper inclusions
 \begin{equation*}
  \Oo^{q;>0}(\hat{S}_d,\EE)\varsubsetneq\Oo^{q;k_1}(\hat{S}_d,\EE)\varsubsetneq\Oo^{q;k_2}(\hat{S}_d,\EE)\varsubsetneq\Oo^{>0}(\hat{S}_d,\EE)\varsubsetneq\Oo^{l_1}(\hat{S}_d,\EE)\varsubsetneq\Oo^{l_2}(\hat{S}_d,\EE).
 \end{equation*}

\end{Rem}

It is also possible to formulate an analogue of Proposition \ref{pr:1}. 
\begin{Prop}
\label{pr:2}
Let $f$ be an entire function and $\hat{f}$ its Taylor expansion, i.e.,
$$
\hat{f}(t)=\sum_{n=0}^{\infty}f_nt^n=\sum_{n=0}^{\infty}\frac{f^{(n)}(0)}{n!}t^n.
$$
Then for every $k>0$
\begin{equation*}
 f\in\Oo^{q;k}(\CC,\EE)\quad\text{if and only if}\quad\hat{f}\in\EE[[t]]_{q;-1/k}.
\end{equation*}
Moreover,
\begin{equation*}
 f\in\Oo^{q;>0}(\CC,\EE)\quad\text{if and only if}\quad\hat{f}\in\EE[[t]]_{q;-\infty}.
\end{equation*}
\end{Prop}

\begin{Def}[{\cite[Definition 14]{L-M}}]
Let $q>1$. A sequence $m$ is said to \emph{preserve $q$-Gevrey asymptotic expansions} if for every $s>0$, $d\in\RR$ and $\hat{u}\in\mathbb{E}[[t]]$ the following statements turn out to be equivalent:
\begin{itemize}
\item[(i)] $\mathcal{B}_{q;s}(\hat{u})\in\mathbb{E}\{t\}$, and this function can be extended on an infinite sector of bisecting direction $d$ with $q$-exponential growth of order $1/s$ on such sector.
\item[(ii)] $\mathcal{B}_{q;s}\mathcal{B}_{m}(\hat{u})\in\mathbb{E}\{t\}$, and this function can be extended on an infinite sector of bisecting direction $d$ with $q$-exponential growth of order $1/s$ on such sector.
\end{itemize}
\end{Def}

\begin{Thm}[{\cite[Theorem 2]{L-M}}]
\label{th:2}
Let $q>1$. A sequence $m=(m(n))_{n\ge0}$ preserves $q$-Gevrey asymptotic expansions if and only if for every $s>0$ and every $\theta\neq 0 \mod 2\pi$, $\mathcal{B}_{m}\left(\sum_{n\ge0}t^n\right)$ and $\mathcal{B}_{m^{-1}}\left(\sum_{n\ge0}t^n\right)$ belong to $\mathbb{C}\{t\}$ and each of them can be extended to an infinite sector of bisecting direction $\theta$ with $q$-exponential growth of order $1/s$. Here, $m^{-1}$ stands for the sequence $(m(n)^{-1})_{n\ge0}$.
\end{Thm}

\begin{Rem}
Observe that every sequence of positive real numbers $m$ which preserves $q$-Gevrey asymptotic expansions is also a sequence preserving summability.
\end{Rem}

Let us now introduce several known examples of sequences preserving Gevrey asymptotic expansions.

\begin{Ex}[{\cite[Proposition 7]{L-M}}]
 If $a>0$ then the sequence $\mathbf{a}=(a^n)_{n\geq 0}$ preserves both Gevrey asymptotic expansions and summability.
\end{Ex}

\begin{Ex}[{\cite[Proposition 8]{L-M}}]
The sequence $m=\left(\frac{(2n)!}{n!^2}\right)_{n\ge0}$ preserves both $q$-Gevrey asymptotic expansions and summability.
\end{Ex}

\begin{Ex}[{\cite[Theorem 4]{L-M}}]
Let $q>1$. The sequence $([n]_{1/q}!)_{n\ge0}$ preserves both $q$-Gevrey asymptotic expansions and summability. 
\end{Ex}

\section{Open set property}\label{sec:3}
In this section we prove that the set of sequences preserving summability (resp. preserving $q$-Gevrey asymptotic expansions) has the following open set property:
\begin{quote}
    if a sequence $\tilde{m}$ belongs to this set then every sequence $m$ close to $\tilde{m}$ also belongs to the same set.
\end{quote}
 It allows us to solve the conjecture posed in \cite{L-M} that the set of sequences preserving $q$-Gevrey asymptotic expansions is strictly contained in the set of sequences preserving summability.

To prove this open set property we introduce a new concept of sequences close to sequences preserving summability (resp. preserving $q$-Gevrey asymptotic expansions).  
\subsection{Sequences close to sequences preserving summability}
 \begin{Thm}
 \label{th:3}
  Let $\tilde{m}=(\tilde{m}(n))_{n\geq0}$ be a sequence preserving summability. We also assume that $m=(m(n))_{n\geq0}$ is a sequence of positive numbers with $m(0)=1$, which is close to the sequence $\tilde{m}$ in the sense that for every $k>0$ there exist $A_k,B_k>0$ satisfying
  \begin{equation}
  \label{eq:eq_2}
   \Big|m(n)-\tilde{m}(n)\Big|\leq\frac{A_kB_k^n}{(n!)^{1/k}}\quad\text{for every}\quad n\in\NN_0.
  \end{equation}
 Then $m$ is also a sequence preserving summability.
 \end{Thm}
 \begin{proof}
  By Theorem \ref{th:1} it is sufficient to show \eqref{eq:eq_1}. To this end we fix $k>0$ and $\theta\neq 0\mod 2\pi$.
  
  In the first part of the proof we will show that
  \begin{equation}\label{eq:eq_2.5}
  \Bo_{m^{-1}}\big(\sum_{n=0}^{\infty}t^n\big)\in\Oo^k(\hat{S}_{\theta}).
  \end{equation}
  Observe that
  \begin{equation}
  \label{eq:eq_3}
   \Bo_{m^{-1}}\big(\sum_{n=0}^{\infty}t^n\big)=\sum_{n=0}^{\infty}m(n)t^n=\sum_{n=0}^{\infty}\big(m(n)-\tilde{m}(n)\big)t^n+\sum_{n=0}^{\infty}\tilde{m}(n)t^n.
  \end{equation}
  Since the sequence $\tilde{m}$ preserves summability,
  by \eqref{eq:eq_1} we see that
  \begin{equation}
  \label{eq:eq_4}
  \sum_{n=0}^{\infty}\tilde{m}(n)t^n=\Bo_{\tilde{m}^{-1}}\big(\sum_{n=0}^{\infty}t^n\big)\in\Oo^k(\hat{S}_{\theta}). 
  \end{equation}
Moreover, by \eqref{eq:eq_2} there exist $A,B>0$ such that
\begin{equation*}
 \Big|\sum_{n=0}^{\infty}\big(m(n)-\tilde{m}(n)\big)t^n\Big|\leq \sum_{n=0}^{\infty}\frac{A_kB_k^n}{(n!)^{1/k}}|t|^n\leq Ae^{B|t|^k}\quad\text{for every}\quad t\in\CC.
\end{equation*}
It means that
\begin{equation}
\label{eq:eq_4.3}
 \sum_{n=0}^{\infty}\big(m(n)-\tilde{m}(n)\big)t^n\in\Oo^k(\CC).
\end{equation}
Therefore by \eqref{eq:eq_3}, \eqref{eq:eq_4} and \eqref{eq:eq_4.3} we deduce that \eqref{eq:eq_2.5} holds

In the second part of the proof we will prove that also
\begin{equation}
\label{eq:eq_4.5}
 \Bo_{m}\big(\sum_{n=0}^{\infty}t^n\big)\in\Oo^k(\hat{S}_{\theta}).
\end{equation}
By Remark \ref{re:1} we see that $\tilde{m}$ is a sequence of order zero. Hence, by Definition~\ref{df:order} and by~\eqref{eq:eq_2} we get
\begin{equation*}
 \Big|m(n)\Big|=\Big|m(n)-\tilde{m}(n)+\tilde{m}(n)\Big|\geq \Big|\tilde{m}(n)\Big|-\Big|m(n)-\tilde{m}(n)\Big|\geq
 a^n-\frac{A_kB_k^n}{(n!)^{1/k}}.
\end{equation*}
Since $\lim\limits_{n\to\infty}\sqrt[n]{A_kB_k^n/(n!)^{1/k}}=0$, there exists $n_0\in\NN_0$ and $\tilde{a}>0$ such that
\begin{equation}
\label{eq:eq_5}
 \Big|m(n)\Big|\geq \tilde{a}^n\quad\text{for every}\quad n\ge n_0.
\end{equation}
We know that $m$ is a sequence of positive numbers with $m(0)=1$. Therefore, decreasing $\tilde{a}>0$ if necessary, we conclude that inequality~\eqref{eq:eq_5} holds for every $n\in\NN_0$. Moreover, by Remark~!\ref{re:1} the sequence $\tilde{m}$ is of order zero. Therefore, using also the assumption~\eqref{eq:eq_2}, we get 
\begin{equation}
\label{eq:eq_8}
 \Big|\frac{1}{m(n)}-\frac{1}{\tilde{m}(n)}\Big|=\Big|\frac{\tilde{m}(n)-m(n)}{\tilde{m}(n)m(n)}\Big|\leq\frac{A_kB_k^n}{\tilde{a}^na^n(n!)^{1/k}}=\frac{A_k\tilde{B}_k^n}{(n!)^{1/k}}\quad\text{for every}\quad n\in\NN_0,
\end{equation}
where $\tilde{B}_k:=\frac{B_k}{\tilde{a}a}$ is a positive constant.

Using the inequality \eqref{eq:eq_8} instead of \eqref{eq:eq_2} and repeating the first part of the proof with $m$ and $\tilde{m}$ replaced by $m^{-1}$ and $\tilde{m}^{-1}$, respectively, we conclude that \eqref{eq:eq_8} holds. 

Combining \eqref{eq:eq_2.5} and \eqref{eq:eq_8} with Theorem \ref{th:1} enables us to obtain the assertion.
\end{proof}

\begin{Rem}
 Let $r(n):=m(n)-\tilde{m}(n)$ for $n\in\NN_0$. Then condition \eqref{eq:eq_2} means that $\sum_{n=0}^{\infty}r(n)t^n\in\CC[[t]]_{-1/k}$ for every $k>0$,
 or, equivalently, that $\sum_{n=0}^{\infty}r(n)t^n\in\CC[[t]]_{-\infty}$. 
\end{Rem}

Using the approach from above we can formulate the following supplement to Theorem~\ref{th:3}:
\begin{Thm}
\label{th:4}
  Let $\tilde{m}=(\tilde{m}(n))_{n\geq0}$ be a sequence preserving summability. We also assume that $m=(m(n))_{n\geq0}$ is a sequence of positive numbers with $m(0)=1$ such that
  there exist $k>0$ satisfying
  \begin{equation}
  \label{eq:eq_8.5}
   \sum_{n=0}^{\infty}r(n)t^n\in\CC[[t]]_{-1/k}\setminus\CC[[t]]_{-\infty},\quad\text{where}\ r(n)=m(n)-\tilde{m}(n)\ \text{for}\ n\in\NN_0.
  \end{equation}
 Then the sequence $m$ does not preserve summability.
 \end{Thm}
 \begin{proof}
Let $w(t):=\sum_{n=0}^{\infty}r(n)t^n$. Because $w(t)\in\CC[[t]]_{-1/k}$, by Proposition!~\ref{pr:1} we see that $w(t)\in\Oo^{k}(\CC)$. From~\eqref{eq:eq_8.5} it follows that $w(t)\not\in\CC[[t]]_{-\infty}$. Hence, there exists $\tilde{k}\in(0,k)$ such that $w(t)\not\in\Oo^{\tilde{k}}(\CC)$. Since $w(t)$ is an entire function, after
decreasing a positive $\tilde{k}$ if necessary, we conclude that $w(t)\not\in\Oo^{\tilde{k}}(\CC\setminus\RR_+)$.

Let $u(t):=\sum_{n=0}^{\infty}m(n)t^n$ and $v(t):=\sum_{n=0}^{\infty}\tilde{m}(n)t^n$. We know that the sequence $\tilde{m}$ preserves summability. Hence, in particular $v(t)\in\Oo^{>0}(\CC\setminus\RR_+)$. On the other hand, $w(t)\not\in\Oo^{\tilde{k}}(\CC\setminus\RR_+)$. Therefore, the function $u(t)=v(t)+w(t)$ does not belong to the space $\Oo^{\tilde{k}}(\hat{S})$, but it means by Theorem \ref{th:1} (see also Remark \ref{re:4}) that the sequence $m$ does not preserve summability.
\end{proof}

\subsection{Sequences close to sequences preserving $q$-Gevrey asymptotic expansions}
Repeating the proofs of Theorems \ref{th:3} and \ref{th:4} with sequences preserving summability replaced by sequences preserving $q$-Gevrey asymptotic expansions, we get
\begin{Thm}
\label{th:5}
  Let $\tilde{m}=(\tilde{m}(n))_{n\geq0}$ be a sequence preserving $q$-Gevrey asymptotic expansions. We also assume that $m=(m(n))_{n\geq0}$ is a sequence of positive numbers with $m(0)=1$ such that for every $s>0$ there exist $A_s,B_s>0$ and $\alpha_s\in\RR$ satisfying
  \begin{equation*}
   \Big|m(n)-\tilde{m}(n)\Big|\leq\frac{A_sB_s^n(n!)^{\alpha_s}}{q^{s\frac{n(n-1)}{2}}}\quad\text{for every}\quad n\in\NN_0.
  \end{equation*}
 Then $m$ is also a sequence preserving $q$-Gevrey asymptotic expansions.
 \end{Thm}

 \begin{Thm}
 \label{th:6}
  Let $\tilde{m}=(\tilde{m}(n))_{n\geq0}$ be a sequence preserving $q$-Gevrey asymptotic expansions. We also assume that $m=(m(n))_{n\geq0}$ is a sequence of positive numbers with $m(0)=1$ such that
  there exist $s>0$ satisfying
  \begin{equation*}
   \sum_{n=0}^{\infty}r(n)t^n\in\CC[[t]]_{q;-s}\setminus\CC[[t]]_{q;-\infty},\quad\text{where}\ r(n)=m(n)-\tilde{m}(n)\ \ \text{for}\ \ n\in\NN_0.
  \end{equation*}
 Then the sequence $m$ does not preserve $q$-Gevrey asymptotic expansions.
 \end{Thm}
 
 Now we are ready to solve the open problem posed in~\cite{L-M} and to describe the relations between the set of sequences
 preserving $q$-Gevrey asymptotic expansions and the set of sequences preserving summability.
 \begin{Thm}
 \label{th:7}
 The set of sequences preserving $q$-Gevrey asymptotic expansions is strictly contained in the set of sequences preserving summability.
 \end{Thm}
 \begin{proof}
  Since $\Oo^{q;>0}(\hat{S}_d,\EE)\subseteq\Oo^{>0}(\hat{S}_d,\EE)$,
  using Theorems \ref{th:1} and \ref{th:2} we conclude that the set of sequences preserving $q$-Gevrey asymptotic expansions is contained in the set of sequences preserving summability.
  
  Let us consider the sequence $m(n)=1+nq^{-\frac{n(n-1)}{2}}$ for $n\in\NN_0$. Since the sequence $\tilde{m}(n)=1^n=\mathbf{1}$ preserves summability, by Theorem \ref{th:3} also the sequence $m=(m(n))_{n\in\NN_0}$ preserves summability. On the other hand, since 
  \begin{equation*}
  \sum_{n=0}^{\infty}nq^{-\frac{n(n-1)}{2}}t^n\in\CC[[t]]_{q;-1}\setminus \CC[[t]]_{q;-\infty},
  \end{equation*}
  by Theorem \ref{th:6} we conclude that the sequence $m$ does not preserve $q$-Gevrey asymptotic expansions. Therefore our inclusion is strict.
 \end{proof}

\section{Special forms of sequences preserving summability}\label{sec:4}

In this chapter we will introduce several classes of sequences preserving summability. For each of them a slightly different approach to showing this fact will be presented.

\subsection{Polynomials}
 Let us consider a complex polynomial $w_p$ of a fixed order $p\in\NN$, satisfying $w_p(0)=1$ and such that $w_p(n)>0$ for every $n\in\NN_0$. It is easy to notice that any such polynomial can be expressed in the form:
\begin{equation}
w_p(n)=a_p(n+1)_p+a_{p-1}(n+1)_{p-1}+\ldots +a_1(n+1)_1+a_0,\label{eq:eq_9}
\end{equation}
where $a_0,\ldots,a_p$ are fixed constant coefficients and $(x)_j$ denotes the Pochhammer symbol, i.e., the increasing factorial
\begin{equation*}
(x)_j=\frac{\Gamma(x+j)}{\Gamma(x)}.
\end{equation*}
It means that
\begin{equation}
	(n+1)_j=(n+1)(n+2)\cdot\ldots\cdot(n+j)\textrm{ for every }j=1,2,\ldots,p.\label{eq:eq_10}
\end{equation}

\begin{Prop}\label{prop:prop3}
Let $m=(m(n))_{n\ge 0}$ be a sequence such that $m(n)=w_p(n)$ for every $n\in\NN_0$. Then for any $k>0$ and $\theta\ne 0\ \mathrm{mod}\ 2\pi$ there exists a disc-sector $\hat{S}_\theta$ such that
\begin{equation}
	\sum_{n=0}^\infty m(n) t^n\in\Oo^k(\hat{S}_\theta).\label{eq:eq_11}
\end{equation}
\end{Prop}
\begin{proof}
    From \eqref{eq:eq_9} and \eqref{eq:eq_10} we conclude that it suffices to prove the theorem for a~sequence $(n+1)_j$ for any fixed $j=1,2,\ldots,p$. Let us notice that, for $|t|<1$, the power series $\sum_{n=0}^\infty (n+1)_jt^n$ is convergent and
	$$
	\sum_{n=0}^\infty (n+1)_jt^n=\frac{d^j}{dt^j}\left(\sum_{n=0}^\infty t^n\right)=\frac{d^j}{dt^j}\left(\frac{1}{1-t}\right)=\frac{ j!}{(1-t)^{j+1}}
	$$
	
	The function  on the right-hand side of the last equality can be analytically extended in any direction except $\theta=0\ \mathrm{mod }\ 2\pi$ and its analytic continuation is bounded at the infinity.
\end{proof}
 
Now, let us notice that any sequence of the form $\frac{1}{w_p(n)}$ can be expressed by means of partial fraction decomposition as
$$
\frac{1}{w_p(n)}=\sum_{q=1}^j\sum_{\nu=1}^{s_q}\frac{C_{q,\nu}}{(n+\alpha_q)^\nu},
$$
where $\alpha_q\in\CC$ for $q=1,\ldots,j$, $C_{q,\nu}\in\CC$ for $\nu=1,\ldots,s_q$ and $q=1,\ldots,j$, and $s_1,\ldots,s_j\in \NN$ such that $s_1+\ldots+s_j=p$. To prove that
\begin{equation}
	\sum_{n=0}^\infty \frac{t^n}{w_p(n)}\in\Oo^k(\hat{S}_\theta)\label{eq:eq_12}
\end{equation}
for any $k>0$ and $\theta\ne 0\ \mathrm{mod}\ 2\pi$, it is enough to come to this conclusion for any sequence of the form $\sum_{n=0}^\infty\frac{t^n}{(n+a)^s}$, where $a\in\CC$ and $s\in\NN$. To arrive at this result, we shall use the Hurwitz-Lerch transcendent and its various properties.

The Hurwitz-Lerch transcendent (see \cite[Section 1.11]{E-M-O-T}) is a function given by the power series
\begin{equation}
	\Phi(t,s,a)=\sum_{n=0}^\infty \frac{t^n}{(n+a)^s}\label{eq:eq_13}
\end{equation}
with $t,s,a\in\CC$, $a\ne 0,-1,-2,\ldots$ and either $|t|<1$ or $|t|=1$ and $\RE s>1$. We shall now focus our attention on the fact that the function given by (\ref{eq:eq_13}) can be analytically extended to $\CC\setminus[1,\infty]$ for fixed values of parameters $s$ and $a$ as defined above.

\begin{Prop}{(see \cite[Lemma 2.1]{G-S})} Consider $a$ such that $\RE a>0$. For $|t|<1$ and $\RE s>0$ or $t=1$ and $\RE s>1$ we have
\begin{equation}
	\Phi(t,s,a)=\frac{1}{\Gamma(s)}\int_0\in\CC[t]^\infty \frac{\zeta^{s-1}\exp(-a\zeta)}{1-t\exp(-\zeta)}d\zeta.\label{eq:eq_14}
\end{equation}
\label{prop:prop4}\end{Prop}

\begin{proof}
	First, let us notice that
	$$
	\frac{1}{1-t\exp(-\zeta)}=\sum_{n\ge 0}t^n\exp(-n\zeta).
	$$
	After substituting the right-hand side of this equality into the integral from~(\ref{eq:eq_14}) we receive
	\begin{multline*}
	\int_0^\infty \zeta^{s-1}\exp(-a\zeta)\left(\sum_{n\ge 0}t^n\exp(-\zeta n)\right)d\zeta=\sum_{n\ge 0}t^n \int_0^\infty \zeta^{s-1}\exp(-\zeta(n+a))\,d\zeta\\=\sum_{n\ge 0}\frac{t^n}{(n+a)^s}\int_0^\infty \zeta^{s-1}\exp(-\zeta)\,d\zeta=\Phi(t,s,a)\Gamma(s),
	\end{multline*}
	which concludes the proof.
\end{proof}

\begin{Prop}{(see \cite[Lemma 2.2]{G-S})} Under previous assumptions about the parameter $a$, and assuming $\RE s>0$, the power series given by~(\ref{eq:eq_13}) can be analytically extended for all $t\in\CC\setminus[1,\infty)$ to an analytic function given by~(\ref{eq:eq_14}).\label{prop:prop5}
\end{Prop}

\begin{Rem}
    Notice that in the framework of the current problem it suffices to consider only $s\in\NN$.
\end{Rem}

The integral representation of the Hurwitz-Lerch transcendent given by~(\ref{eq:eq_14}) can also be made valid for a fixed $a$ such that $\RE a<0$. First of all, let us notice that all $a$'s that we consider are such that $-a$'s are zeros of the polynomial $w_p(n)$, about which we assumed initially that $w_p(0)=1$ and $w_p(n)\ne 0$ for every $n\in\NN$. Hence, $-a\not\in\NN_0$. Moreover, for any $a$ satisfying $\RE a<0$ there exists $n_0=n_0(a)$ such that $n_0+\RE a>0$. We can then write
$$
\sum_{n\ge 0}\frac{t^n}{(n+a)^s}=\sum_{n=0}^{n_0-1}\frac{t^n}{(n+a)^s}+\sum_{n\ge n_0}\frac{t^n}{(n+a)^s}=\sum_{n=0}^{n_0-1}\frac{t^n}{(n+a)^s}+t^{n_0}\Phi(t,s,n_0+a),
$$
and $\Phi(t,s,n_0+a)$ can be expressed using~(\ref{eq:eq_14}).

\begin{Prop}\label{prop:prop6}
    For any direction $\theta\neq 0\ \mathrm{mod}\ 2\pi$ the function $\Phi(t,s,a)$ is of less than exponential growth. 
\end{Prop}

\begin{proof}
    Let us fix a direction $\theta\neq 0\ \mathrm{mod}\ 2\pi$ and let $t=re^{i\theta}$. Then
	$$
	\Phi(t,s,a)=\frac{1}{\Gamma(s)}\int_0^{\infty}\frac{\zeta^{s-1}e^{-a\zeta}}{1-re^{-\zeta+i\theta}}\,d\zeta.
	$$
	It is enough to show that $\Phi(re^{i\theta},s,a)$ can be bounded from above by a positive constant the same for all $r>0$.
	
	Notice that
	$$|1-re^{-t+i\theta}|\ge\inf_{\tau>0,\,\rho>0}|1-\rho e^{-\tau+i\theta}|=\mathrm{dist}\{1,L_{\theta}\},$$
	with $\mathrm{dist}\{1,L_{\theta}\}$ denoting the distance on a complex plane between $1$ and the half-line starting at the origin and going in direction $\theta$. Since $$\mathrm{dist}\{1,L_{\theta}\}=\inf_{x>0}|1-xe^{i\theta}|=C(\theta)$$
	with
        \begin{equation}\label{eq:C}
	C(\theta)=\begin{cases}
		|\sin\theta|&\textrm{ for }\theta\in\left(-\frac{\pi}{2},\ \frac{\pi}{2}\right)\setminus\{0\}\textrm{ mod }2\pi\\
		1&\textrm{ otherwise}
	\end{cases},
	\end{equation}
	we receive
	$$
	\left|\int_0^{\infty}\frac{\zeta^{s-1}e^{-a\zeta}}{1-re^{-\zeta+i\theta}}\,d\zeta\right|\le\frac{1}{C(\theta)}\int_0^\infty \zeta^{s-1}e^{-\RE a\zeta}\,d\zeta.
	$$
	After using substitution $\RE a\zeta=w$, we obtain
	$$
	\left|\int_0^{\infty}\frac{\zeta^{s-1}e^{-a\zeta}}{1-re^{-\zeta+i\theta}}\,d\zeta\right|\le\frac{1}{C(\theta)}\frac{\Gamma(s)}{(\RE a)^s}.
	$$
	
	Hence, for any $r>0$ we receive
	$$
	|\Phi(re^{i\theta},s,a)|\le [C(\theta)(\RE a)^s]^{-1}.
	$$
\end{proof}

Using the reasoning above we can prove that a sequence $m(n)$ given by a polynomial is a sequence preserving summability.

\begin{Thm}\label{th:th_8}
    Let $w_p$ be a polynomial with complex coefficients of order $p\in\NN$ satisfying $w_p(n)>0$ for every $n\in\NN_0$ and $w_p(0)=1$. Then, the sequence $m=(m(n))_{n\ge 0}$ defined as $m(n)=w_p(n)$ for every $n\in\NN_0$ preserves summability.
\end{Thm}

\begin{proof}
    As stated in Theorem \ref{th:1}, it is enough to prove that for every $k>0$ and for every $\theta\neq 0\mod 2\pi$ there exists a disc-sector $\hat{S}_{\theta}$ such that $\Bo_{m}\big(\sum_{n=0}^{\infty}t^n\big)\in\Oo^k(\hat{S}_{\theta})$ and $\Bo_{m^{-1}}\big(\sum_{n=0}^{\infty}t^n\big)\in\Oo^k(\hat{S}_{\theta})$. The first of these facts has been proven in Proposition \ref{prop:prop3}. Therefore, it remains to show, for any fixed $k>0$ and $\theta\neq 0\mod\ 2\pi$, that
    $$
    \sum_{n=0}^\infty\frac{t^n}{w_p(n)}\in\Oo^k(\hat{S}_{\theta}).
    $$
    for some disc-sector $\hat{S}_{\theta}$.

    Since by means of partial fraction decomposition we have
    $$
    \frac{1}{w_p(n)}=\sum_{q=1}^j\sum_{\nu=1}^{s_q}\frac{C_{q,\nu}}{(n+\alpha_q)^\nu},
    $$
    where $w_p(\alpha_q)=0$ for $q=1,\ldots,j$, the conclusion follows directly from Propositions~\ref{prop:prop4}, \ref{prop:prop5} and \ref{prop:prop6}.
\end{proof}

Since the sets of sequences preserving summability forms a group with multiplication (see Remark \ref{re:groups}),
we may extend polynomials given in Theorem \ref{th:th_8} to rational functions. Namely we conclude that
\begin{Cor}
Assume that $w_1,w_2\in\CC[t]$ are polynomials satisfying $w_1(n)/w_2(n)>0$ for every $n\in\NN_0$ and $w_1(0)=w_2(0)$. Then the sequence $m=(m(n))_{n\geq 0}$ defined as $m(n)=w_1(n)/w_2(n)$ for every $n\in\NN_0$ preserves summability.
\end{Cor}


\subsection{Sums of powers}
We show a similar result for sequences being sums of powers.
\begin{Thm}\label{th:9}
Let $m(n)={a_1}^n+{a_2}^n+\ldots+{a_k}^n$ for $n\in\NN_0$, where ${a_1}>{a_2}>\ldots>{a_k}>0$ and $k\in\NN$. Then the sequence $m=(m(n))_{n\geq 0}$ preserves summability.
\end{Thm}
\begin{proof}
Consider the power series 
\begin{multline*}
\sum_{n=0}^\infty m(n) t^n=\sum_{n=0}^\infty ({a_1}^n+{a_2}^n+\ldots+{a_k}^n) t^n=\sum_{n=0}^\infty ({a_1t})^n+\sum_{n=0}^\infty ({a_2t})^n+\ldots+\sum_{n=0}^\infty ({a_kt})^n=\\
=\frac{1}{1-a_1t}+\frac{1}{1-a_2t}+\ldots+\frac{1}{1-a_kt},
\end{multline*}
which is convergent for $|t|<\frac{1}{a_1}$ and $\sum_{n=0}^\infty m(n) t^n\in\Oo(\CC\setminus\RR_+).$

Observe also that for $t\in\CC\setminus\RR_+$ such that $\arg t= \theta\neq 0\mod 2\pi$, we have
$$
\bigg|\sum_{n=0}^\infty m(n) t^n\bigg|\leq \frac{k}{\mathrm{dist}\{1,L_{\theta}\}}=\frac{k}{|\sin\theta|},
$$
where $\mathrm{dist}\{1,L_{\theta}\}$ once again denotes the distance on a complex plane between $1$ and the half-line $L_{\theta}$ starting at the origin and going in direction $\theta$, and $C(\theta)$ is defined by \eqref{eq:C}. Hence, $\sum_{n=0}^\infty m(n) t^n\in\Oo^{>0}(\CC\setminus\RR_+).$

Now, let $f(t)=\sum_{n=0}^\infty \frac{t^n}{m(n)}\in\Oo(D_r)$, where a radius of convergence is given by $r=\lim_{n\longrightarrow\infty}\sqrt[n]{m(n)}=a_1$. Then
$$
\sum_{n=0}^\infty \frac{(a_1t)^n}{m(n)}+\sum_{n=0}^\infty \frac{(a_2t)^n}{m(n)}+\ldots+\sum_{n=0}^\infty \frac{(a_kt)^n}{m(n)}=\sum_{n=0}^\infty t^n=\frac{1}{1-t} ,
$$
hence 
$$
f(a_1t)+f(a_2t)+\ldots+f(a_kt)=\frac{1}{1-t}.
$$
After the substitution $z=a_1t$ is applied, the last equality will take the form
$$
f(z)+f\bigg(\frac{a_2}{a_1}z\bigg)+f\bigg(\frac{a_3}{a_1}z\bigg)+\ldots+f\bigg(\frac{a_k}{a_1}z\bigg)=\frac{1}{1-\big(\frac{z}{a_1}\big)}.
$$
Then, putting $\tau_j=\frac{a_j}{a_1}$ for $j=2,3,\ldots,k$ such that $1>\tau_2>\ldots>\tau_k$, we obtain
$$
f(z)+f(\tau_2z)+f(\tau_3z)+\ldots+f(\tau_kz)=\frac{1}{1-\big(\frac{z}{a_1}\big)},
$$
hence
\begin{equation}
f(z)=-f(\tau_2z)-f(\tau_3z)-\ldots-f(\tau_kz)+\frac{1}{1-\big(\frac{z}{a_1}\big)}.\label{eq:eq_16}
\end{equation}
Repeated use of the formula given above leads us to a conclusion that the function $f(z)\in\Oo\big(D_{r}\cap(\CC\setminus\RR_+)\big)$ can be analytically extended to $\CC\setminus\RR_+$.

Let $z\in\CC\setminus\RR_+$ and $\arg z=\theta\neq 0\mod 2\pi$. Then there exists the smallest $N\in\NN$ such that $\tau_2^Nz\in D_{r/2}$.
Now we proceed with the equation \eqref{eq:eq_16} as follows:
first we replace $z$ by $\tau_2 z$ twice and get
$$
f(z)=\sum_{i=2}^k\sum_{j=2}^kf(\tau_i\tau_j z)+\frac{1}{1-(\frac{z}{a_1})}-\sum_{i=2}^k\frac{1}{1-\big(\frac{\tau_iz}{a_1}\big)},
$$
then we iterate until the Nth step, in which we obtain
\begin{multline*}
f(z)=\sum_{j_1=2}^k\ldots\sum_{j_N=2}^kf(\tau_{j_1}\ldots\tau_{j _N}z)+\frac{1}{1-\big(\frac{z}{a_1}\big)}-\sum_{i=2}^k\frac{1}{1-\big(\frac{\tau_iz}{a_1}\big)}+\ldots\pm\\\pm \sum_{j_1=2}^k\ldots\sum_{j_{N-1}=2}^k\frac{1}{1-\frac{\tau_{j_1}\ldots\tau_{j_{N-1}}z}{a_1}}.
\end{multline*}
Since $\tau_2^N<\frac{r}{2|z|}\leq\tau_2^{N-1}$, then 
\begin{equation}
N-1\leq\log_{\frac{1}{\tau_2}}\Big(\frac{2|z|}{r}\Big)<N.\label{eq:eq_17}
\end{equation}
Let us now notice that
\begin{multline*}
\big|f(z)\big|\leq (k-1)^N\cdot B+ \bigg(1+(k-1)+(k-1)^2+\ldots+(k-1)^{N-1}\bigg)\cdot\frac{1}{\mathrm{dist}\{1,L_{\theta}\}}=\\
=(k-1)^N\cdot B+\frac{(k-1)^N-1}{k-2}\cdot \frac{1}{C(\theta)}\leq(k-1)^N\bigg(B+\frac{1}{C(\theta)}\bigg),
\end{multline*}
where $B=\sup_{|x|\leq r/2}|f(x)|<\infty$.
From \eqref{eq:eq_17} we receive 
$$
N\leq\log_{\frac{1}{\tau_2}}\Big(\frac{2|z|}{r}\Big)+1\ \textrm{ if and only if }\ N\leq\log_{\frac{1}{\tau_2}}\Big(\frac{2|z|}{\tau_2r}\Big)\Longleftrightarrow N\leq\frac{\log_{k-1}\Big(\frac{2|z|}{\tau_2r}\Big)}{\log_{k-1}\Big(\frac{1}{\tau_2}\Big)}.
$$
Using the last inequality we can write
\begin{multline*}
\big|f(z)\big|\leq(k-1)^N\cdot\bigg(B+\frac{1}{C(\theta)}\bigg)\leq(k-1)^{\frac{\log_{k-1}(\frac{2|z|}{\tau_2r})}{\log_{k-1}(\frac{1}{\tau_2})}}\cdot\bigg(B+\frac{1}{C(\theta)}\bigg)\leq \\
\leq\Big(\frac{2}{\tau_2 r}\Big)^{\frac{1}{\log_{k-1}(\frac{1}{\tau_2})}}\cdot|z|^{\frac{1}{\log_{k-1}(\frac{1}{\tau_2})}}\cdot\bigg(B+\frac{1}{C(\theta)}\bigg)\leq
A|z|^{\log_{\frac{1}{\tau_2}}(k-1)},
\end{multline*}
where $A=\Big(\frac{2}{\tau_2 r}\Big)^{\log_{\frac{1}{\tau_2}}(k-1)}\bigg(B+\frac{1}{C(\theta)}\bigg)$. 
Thus, $f(z)\in\Oo^{>0}(\CC\setminus\RR_+)$ and 
by Theorem~\ref{th:1} we conclude that the sequence $(m(n))_{n\geq 0}$ preserves summability. 
\end{proof}

\begin{Ex}
Let $m(n)=1+2^n$ for $n\in\NN_0$. We show directly that the sequence $(m(n))_{n\geq 0}$ preserves summability.  We easily see that $$\sum_{n=0}^\infty m(n) t^n= \sum_{n=0}^\infty (1+2^n) t^n\in\Oo(\CC\setminus\RR_+).$$ Observe also that, for $|t|<\frac{1}{2}$ the power series
$$\sum_{n=0}^\infty m(n) t^n= \sum_{n=0}^\infty (1+2^n) t^n=\sum_{n=0}^\infty t^n+\sum_{n=0}^\infty (2t)^n=\frac{1}{1-t}+\frac{1}{1-2t}$$ is convergent.

Let us now consider a function $f(t)=\sum_{n=0}^\infty \frac{t^n}{m(n)}\in\Oo(D_r)$ with a radius of convergence $r=\lim_{n\longrightarrow\infty}\sqrt[n]{m(n)}=2$. Then
$$
\sum_{n=0}^\infty \frac{t^n}{m(n)}+\sum_{n=0}^\infty \frac{(2t)^n}{m(n)}=\sum_{n=0}^\infty t^n=\frac{1}{1-t}.
$$
Hence, 
$$
f(t)+f(2t)=\frac{1}{1-t}\ \textrm{ if and only if }\  f(2t)=\frac{1}{1-t}-f(t)\in\Oo(D_r). 
$$
After replacing $2t$ by $z$ in the last equation, we obtain
$$
f(z)=\frac{1}{1-\frac{z}{2}}-f\bigg(\frac{z}{2}\bigg),\,\,\,\mathrm{for}\,\, |z|<4\textrm{ and } z\in\Oo(\CC\setminus\RR_+).
$$
Thus the function $f(z)\in\Oo\big(D_{r}\cap(\CC\setminus\RR_+)\big)$ can be analytically extended to $\CC\setminus\RR_+$.

We now proceed to show that $f\in\Oo^{>0}(\CC\setminus\RR_+)$. First let us observe that
\begin{align*}
f(2t)&=\frac{1}{1-t}-f(t),\\\\
f(4t)&=\frac{1}{1-2t}-f(2t)=\frac{1}{1-2t}-\left(\frac{1}{1-t}-f(t)\right)=\frac{1}{1-2t}-\frac{1}{1-t}+f(t),\\
&\qquad\vdots\\
f(2^Nt)&=\frac{1}{1-2^{N-1}t}-\frac{1}{1-2^{N-2}t}+\ldots-\frac{1}{1-t}+f(t),\,\, N\in\NN.
\end{align*}

Let $\tau\in\CC\setminus\RR_+$ and $\arg \tau=\theta\neq 0\mod 2\pi$. Then there exists the smallest $N\in\NN$ such that $\frac{\tau}{2^N}\in D_{\frac{1}{2}r}$.
Therefore, the equation 
$$
f(2^Nt)=\frac{1}{1-2^{N-1}t}-\frac{1}{1-2^{N-2}t}+\ldots-\frac{1}{1-t}+f(t)
$$ 
can be rewritten as
$$
f(\tau)=\sum_{n=1}^N(-1)^{n-1}\frac{1}{1-\frac{\tau}{2^n}}+f\Big(\frac{\tau}{2^N}\Big).
$$ 
Since $\frac{1}{2}\leq\frac{|\tau|}{2^N}< 1$, we obtain 
\begin{equation}
N-1\leq\log_2|\tau|< N.\label{eq:eq_15}
\end{equation}

Now let us notice that 
$$
|f(\tau)|\leq N\cdot\frac{1}{\mathrm{dist}\{1,L_{\theta}\}} + B=\frac{N}{C(\theta)}+B, 
$$ 
where $B=\sup_{|z|<1}|f(z)|$ and $\mathrm{dist}\{1,L_{\theta}\}=C(\theta)$ with $C(\theta)$ given by \eqref{eq:C}.
Thus~\eqref{eq:eq_15} leads to 
$$|f(\tau)|\leq \big(\log_2|\tau|+1\big)\cdot\frac{1}{|\sin \theta|}+ C.$$
We arrive at the desired conclusion, i.e., $f(\tau)\in\Oo^{>0}(\CC\setminus\RR_+)$.
\end{Ex}

\subsection{Solutions of linear difference equations with constant coefficients}
Our next goal is to generalize the previous results to sequences that are solutions of homogeneous linear difference equations with constant coefficients. Namely, we will study sequences $(m(n))_{n\geq 0}$ of positive numbers given in the form
\begin{equation}\label{eq:m_n}
m(n)=w_1(n)a_1^n+\ldots+w_l(n)a_l^n\quad\textrm{for}\quad n\in\NN_0,
\end{equation}
where $a_1,\ldots,a_l$ are positive numbers and $w_1,\ldots,w_l$ are given polynomials with complex coefficients. Observe that for $a_1=\ldots=a_l=1$ we get polynomials considered in Theorem \ref{th:th_8}, and for $w_1(n)\equiv\ldots\equiv w_l(n)\equiv 1$ we receive sums of powers studied in Theorem \ref{th:9}.

In order to obtain a general result for sequences defined by \eqref{eq:m_n}, first we need to prove several technical lemmas.
\begin{Lem}\label{lemat 1}
Let $w(n)=a_p(n+1)^p+\ldots+a_1(n+1)+a_0$ be a polynomial of order $p\in\NN$ with complex coefficients and let $A=|a_0|+|a_1|+\ldots+|a_p|$. Then for every $k\in\NN_0$ the function $f(z)=\sum_{n=0}^{\infty}w^k(n)z^n$ is analytic on the complex disc $D_1$ and can be analytically continued to $\CC\setminus[1,\infty)$. Moreover, this continuation is estimated by
\begin{equation*}
\sup_{z\in L_{\theta}}|f(z)|\leq \frac{A^k\big((kp)!\big)^2}{C(\theta)^{kp+1}}\quad\textrm{for every}\quad\theta\neq 0\mod 2\pi,
\end{equation*}
where $L_{\theta}=\{z\in\CC: \arg z=\theta\}$ and $C(\theta)=\begin{cases}
		|\sin\theta|&\textrm{ for }\theta\in\left(-\frac{\pi}{2},\ \frac{\pi}{2}\right)\setminus\{0\}\textrm{ mod }2\pi\\
		1&\textrm{ otherwise}
	\end{cases}.$
\end{Lem}
\begin{proof}
First, we consider the series $\sum_{n=0}^{\infty}(n+1)_pz^n$ 
. Observe that for $z\in D_1$
\begin{equation}\label{lem1 (1)}
    \sum_{n=0}^{\infty}(n+1)_pz^n=\frac{d^p}{dz^p}\bigg(\sum_{n=0}^{\infty}z^n\bigg)=\frac{d^p}{dz^p}\bigg(\frac{1}{1-z}\bigg)=\frac{p!}{(1-z)^{p+1}},
\end{equation}
so $f_1(z)=\sum_{n=0}^{\infty}(n+1)_pz^n$ can be analytically continued to $\CC\setminus[1,\infty)$ and $$\sup_{z\in L_{\theta}}|f_1(z)|\leq\frac{p!}{C(\theta)^{p+1}}.$$

Next, we assume that $f_2(z)=\sum_{n=0}^{\infty}(n+1)^{p}z^n$ and we use the relation between powers and rising factorials (see \cite[formula (6.12)]{G-K-P})
\begin{equation}\label{lem1 (2)}
    a^p=\sum_{m=0}^p \bigg\{{p\atop m}\bigg\}(-1)^{p-m}(a)_m,
    \end{equation}
where $\Big\{{p\atop m}\Big\}$ are Stirling numbers of the second kind.
We also use the following properties of Stirling numbers (see \cite[(6.6) and (6.9)]{G-K-P})
\begin{equation}\label{lem1 (3)}
    \Big\{{p\atop m}\Big\}\leq \Big[{p\atop m}\Big]\quad\mathrm{and}\quad\sum_{m=0}^{p}\Big[{p\atop m}\Big]=p!,
\end{equation}
where $\Big[{p\atop m}\Big]$ are Stirling numbers of the first kind.
By (\ref{lem1 (1)}) and (\ref{lem1 (2)}) we get
\begin{equation}\label{lem1 (4)}
    \sum_{n=0}^{\infty}(n+1)^pz^n=\sum_{m=0}^p \bigg\{{p\atop m}\bigg\}(-1)^{p-m}\sum_{n=0}^{\infty}(n+1)_mz^n=\sum_{m=0}^p \bigg\{{p\atop m}\bigg\}(-1)^{p-m}\frac{m!}{(1-z)^{m+1}}.
\end{equation}
Hence, using (\ref{lem1 (3)}) we estimate $f_2(z)=\sum_{n=0}^{\infty}(n+1)^pz^n\in\Oo(\CC\setminus[1,\infty))$ as follows
\begin{multline*}  
\sup_{z\in L_{\theta}}|f_2(z)|\leq \sum_{m=0}^p \bigg\{{p\atop m}\bigg\}\frac{m!}{C(\theta)^{m+1}}\leq \frac{p!}{C(\theta)^{p+1}}\sum_{m=0}^p \bigg\{{p\atop m}\bigg\}\leq\\\leq\frac{p!}{C(\theta)^{p+1}}\sum_{m=0}^p \bigg[{p\atop m}\bigg]=\frac{(p!)^2}{C(\theta)^{p+1}}.
\end{multline*}

Next, we take $f_3(z)=\sum_{n=0}^{\infty}w(n)z^n$, where $w(n)=a_p(n+1)^p+\ldots+a_1(n+1)+a_0.$
From (\ref{lem1 (4)}) it follows that
\begin{equation}\label{lem1 (5)}
f_3(z)=\sum_{n=0}^{\infty}w(n)z^n=\sum_{l=0}^p a_l\sum_{m=0}^l\bigg\{{l\atop m}\bigg\}(-1)^{l-m}\frac{m!}{(1-z)^{m+1}}.
\end{equation}
It means that $f_3(z)$ can be analytically continued to $\CC\setminus[1,\infty)$ and 
\begin{equation}\label{lem1 (6)}
\sup_{z\in L_{\theta}}|f_3(z)|\leq\frac{A(p!)^2}{C(\theta)^{p+1}},
\end{equation}

Finally, let us observe that $w^k(n)=\big[a_p(n+1)^p+\ldots+a_1(n+1)+a_0\big]^k$ is a~polynomial of order $kp$ of the form $w^k(n)=b_{kp}(n+1)^{kp}+\ldots+b_1(n+1)+b_0$, where $|b_0|+|b_1|+\ldots|b_{kp}|=\big(|a_0|+|a_1|+\ldots+|a_p|\big)^k=A^k$.

Hence, after replacing $w(n)$ by $w^k(n)$ in (\ref{lem1 (5)}) we see that
$$f(z)=\sum_{n=0}^{\infty}w^k(n)z^n=\sum_{l=0}^{kp}b_l\sum_{m=0}^ l\bigg\{{l\atop m}\bigg\}(-1)^{l-m}\frac{m!}{(1-z)^{m+1}}.$$ Therefore $f(z)$ can be analytically continued to $\CC\setminus[1,\infty)$ and from~(\ref{lem1 (6)}) it follows that
$$\sup_{z\in L_{\theta}}|f(z)|\leq\frac{A^k\big((kp)!\big)^2}{C(\theta)^{kp+1}}.$$
\end{proof}

\begin{Lem}\label{lemat 2}
    Assume that $f(z)=\sum_{n=0}^{\infty}u_nz^n\in\Oo(D_r)$ for some $r>0$, $f(z)$ can be analytically continued to $\CC\setminus[r,\infty)$ and for every $\theta\neq 0\mod 2\pi$ there exists a non-decreasing function $p_{\theta}:[0,\infty)\longrightarrow [0,\infty)$ such that 
    $$|f(z)|\leq p_\theta(|z|)\quad for \quad z\in L_\theta.$$
    Then for every $a\in\CC$ with $\mathrm{Re}\,a>1$ the function $g(z)=\sum_{n=0}^{\infty}\frac{u_n}{(n+a)}z^n$ is analytic on $D_r$ and can be analytically continued to $\CC\setminus[r,\infty)$. Moreover, for every $\theta\neq 0\mod 2\pi$ $$|g(z)|\leq p_\theta(|z|)\quad for \quad z\in L_\theta.$$
\end{Lem}
\begin{proof}
First, let us observe that by the Cauchy-Hadamard Theorem the power series $g(z)$ has the same radius of convergence as $f(z)$, so $g(z)\in\Oo(D_r)$. Additionally for $0<|z|<r$ we get
\begin{multline*}g(z)=\sum_{n=0}^{\infty}\frac{u_n}{(n+a)}z^n=z^{-a}\sum_{n=0}^{\infty}\frac{u_n}{(n+a)}z^{n+a}=\\=z^{-a}\int_0^z\sum_{n=0}^{\infty}u_n\zeta^{n+a-1}d\zeta=z^{-a}\int_0^z\zeta^{a-1}f(\zeta)d\zeta.
\end{multline*}
Since $f(z)$ can be analytically continued to $\CC\setminus[r,\infty)$, by the above integral representation of $g(z)$ we conclude that also $g(z)$ can be analytically continued to $\CC\setminus[r,\infty)$. We may also estimate 
\begin{multline*}
    |g(z)|\leq \int_0^{|z|}\bigg|\frac{(re^{i\theta})^a}{\big(|z|e^{i\theta}\big)^a}\bigg|\cdot\big|re^{i\theta}\big|^{-1}\cdot\big|f(re^{i\theta})\big|dr\leq \frac{1}{|z|^a}\int_0^{|z|}r^{a-1}p_\theta(r)dr\leq\\\leq\frac{|z|^{a-1}p_\theta(|z|)}{|z|^a}\int_0^{|z|}dr\leq p_\theta(|z|)
\end{multline*}
for $z\in L_\theta$.
\end{proof}

\begin{Lem}\label{lemat 3}
Let $r(n)=\frac{w_1(n)}{w_2(n)}=w(n)+q(n)$ be a general rational function, where $w_2(n)\neq 0$ for $n\in\NN_0$, $w(n)=a_p(n+1)^p+\ldots+a_1(n+1)+a_0$ and $q(n)=\sum_{q=1}^m\frac{C_q}{(n+\tilde a_p)^{{\nu}_q}}$, where ${\nu}_q\in\NN$ and $C_q, \tilde a_p\in\CC$. Then for every $k\in\NN$ function $f(z)=\sum_{n=0}^\infty r^k(n)z^n$ is analytic on $D_1$ and can be analytically continued to $\CC\setminus[1,\infty)$. Moreover, for every $\theta\neq 0\mod 2\pi$ there exist constants $C<\infty$ and $n_0\in\NN_0$ dependent only on $r(n)$ such that $$|f(z)|\leq\frac{C^k \big((kp)!\big)^2}{C(\theta)^{pk+1}}\left(|z|^{n_0}+1\right)\quad for \ z\in L_\theta.$$
\end{Lem}
\begin{proof} We shall divide our proof into three steps.

\emph{Step 1.} Let $$\tilde f(z)=\sum_{n=0}^{\infty}w^{l_0}(n)\prod_{j=1}^l \frac{\tilde C_j}{(n+\tilde a_j)}z^n$$
for some $l_0,l\in\NN_0$ and $\tilde a_j\in\CC$ with $\mathrm{Re}\, \tilde a_j>1$ $(j=1,2,\ldots,l)$.
Observe also that by Lemmas~\ref{lemat 1} and~\ref{lemat 2}, for$\tilde f(z)\in\Oo(\CC\setminus[1,\infty))$ and moreover, for $g(z)=\sum_{n=0}^{\infty}w^{l_0}(n)z^n$ we have
$$\sup_{z\in L_\theta}|\tilde f(z)|\leq \tilde C_1\cdot \tilde C_2\cdot\ldots\cdot \tilde C_l\,\,\sup_{z\in L_\theta}|g(z)|.$$ Hence, by Lemma \ref{lemat 1}
\begin{equation}\label{lem3 (7)}
|\tilde f(z)|\leq \frac{\tilde C^l\cdot A^{l_0}\big((l_0p)!\big)^2}{C(\theta)^{l_0p+1}},\quad\mathrm{for}\ z\in L_\theta,    
\end{equation}
where $A=|a_0|+|a_1|+\ldots+|a_p|$ and $\tilde C=\max_{j=1,2,\ldots,l}|\tilde C_j|$.
\medskip\par
\emph{Step 2.} Since $$\left(w(n)+\sum_{q=1}^m\frac{C_q}{(n+\tilde a_q)^{{\nu}_q}}\right)^k=\sum_{l_0+l_1+\ldots+l_m=k}\frac{k!}{l_0!\cdot\ldots\cdot l_m!}w^{l_0}(n)\prod_{q=1}^m\left(\frac{C_q}{(n+\tilde a_q)^{{\nu}_q}}\right)^l,$$
using (\ref{lem3 (7)}) we conclude that
\begin{multline*}
\big|f(z)\big|=\Bigg|\sum_{n=0}^\infty\bigg(w(n)+\sum_{q=1}^m\frac{C_q}{(n+\tilde a_q)^{{\nu}_q}}\bigg)^kz^n\Bigg|\leq \\ \leq \sum_{l_0+l_1+\ldots+l_m=k}\frac{k!}{l_0!\cdot\ldots\cdot l_m!}\Bigg|\sum_{n=0}^\infty w^{l_0}(n)\cdot \prod_{q=1}^m\bigg(\frac{C_q}{(n+\tilde a_q)^{{\nu}_q}}\bigg)^{l_q}z^n\Bigg|\leq \\ \leq\sum_{l_0+l_1+\ldots+l_m=k}\frac{k!}{l_0!\cdot\ldots\cdot l_m!} \frac{A^{l_0}|C_1|^{l_1}\cdot\ldots\cdot |C_m|^{l_m}\big((l_0p)!\big)^2}{C(\theta)^{l_0p+1}}\leq \\ \leq \frac{\big((kp)!\big)^2}{C(\theta)^{kp+1}}\Big(A+|C_1|+\ldots+|C_m|\Big)^k.
\end{multline*}
Hence, $$\big|f(z)\big|\leq \frac{C^k ((kp)!)^2}{C(\theta)^{kp+1}}\quad\mathrm{for}\quad z\in L_\theta,\quad\mathrm{where}\quad C=A+|C_1|+\ldots+|C_m|,$$
under condition that $\mathrm{Re}\,\tilde a_q>1$ for $q=1,\ldots,m$.
\medskip\par
\emph{Step 3.}
In the general case we may find $n_0\in\NN$ such that 
\begin{equation}
\label{lem3 (8)}\mathrm{Re}\,(n_0+\tilde a_q)>1\quad\mathrm{for}\quad q=1,\ldots,m.
\end{equation}
Then $$f(z)=\sum_{n=0}^\infty r^k(n)z^n=\sum_{n=0}^{n_0-1}r^k(n)z^n+\sum_{n=n_0}^\infty r^k(n)z^n.$$
To estimate the first sum observe that
$$\bigg|\sum_{n=0}^{n_0-1}r^k(n)z^n\bigg|\leq \tilde A^k\Big(|z|^{n_0-1}+1\Big),\quad\mathrm{where}\quad \tilde A=|r(0)|+\ldots+|r(n_0-1)|.$$
We also have
$$\bigg|\sum_{n=n_0}^\infty r^k(n)z^n\bigg|=|z|^{n_0}\cdot\bigg|\sum_{n=0}^\infty r(n+n_0)^pz^n\bigg|.$$
Since $$r(n+n_0)=w(n+n_0)+\sum_{q=1}^m\frac{C_q}{(n+n_0+\tilde a_q)^{{\nu}_q}}$$ and (\ref{lem3 (8)}) holds, using the previous parts of the proof we conclude that
$$\bigg|\sum_{n=n_0}^\infty r^k(n)z^n\bigg|\leq|z|^{n_0}\cdot\frac{\bar C^k\big((kp)!\big)^2}{C(\theta)^{kp+1}}\quad\mathrm{for}\quad z\in L_\theta,$$
where $C=A+|C_1|+\ldots+|C_m|.$

Hence, we finally arrive at
$$|f(z)|\leq\frac{C^k\big((kp)!\big)^2}{C(\theta)^{kp+1}} \Big(|z|^{n_0}+1\Big)\quad\mathrm{for}\quad z\in L_\theta,\quad\mathrm{where}\quad C=A+|C_1|+\ldots+|C_m|+\tilde A.$$
\end{proof}

Now we are ready to prove that sequences $(m(n))_{n\geq 0}$ defined by \eqref{eq:m_n} preserve summability. 

\begin{Thm}\label{th:10}
Assume that $m=(m(n))_{n\ge 0}$ is a sequence of positive numbers defined by the formula $m(n)=w_1(n)a_1^n+\ldots+w_l(n)a_l^n$ for every $n\in\NN_0$ , where $0<a_1<\ldots<a_l$ and $w_1,\ldots,w_l$ are given polynomials with complex coefficients, and $w_l(n)\neq0$ for every $n\in\NN_0$. Then the sequence $m=(m(n))_{n\in\NN_0}$ preserves summability.
\end{Thm}
\begin{proof}
First, let us observe that for $|z|<\frac{1}{a_l}$
\begin{multline*}
\sum_{n=0}^\infty m(n)z^n=w_1\Big(z\frac{d}{dz}\Big)\sum_{n=0}^\infty(a_1z)^n+\ldots+w_l\Big(z\frac{d}{dz}\Big)\sum_{n=0}^\infty(a_lz)^n=\\= w_1\Big(z\frac{d}{dz}\Big)\frac{1}{1-a_1z}+\ldots+ w_l\Big(z\frac{d}{dz}\Big)\frac{1}{1-a_lz}.
\end{multline*}
Hence, this function can be analytically continued to $\CC\setminus\big[\frac{1}{a_l},\infty)$ and its analytic continuation belongs to the space $\Oo^{>0}\Big(\CC\setminus\big[\frac{1}{a_l},\infty)\Big)$. 

Now, we consider the function $f(z)=\sum_{n=0}^\infty \frac{z^n}{m(n)}$, which is holomorphic for $|z|<a_l$. Observe that $$w_1\Big(z\frac{d}{dz}\Big)f(a_1z)+\ldots+w_l\Big(z\frac{d}{dz}\Big)f(a_lz)=\sum_{n=0}^\infty \frac{w_1(n)a_1^n+\ldots+w_l(n)a_l^n}{m(n)}z^n=\frac{1}{1-z},$$
and so we can see that $$w_l\Big(z\frac{d}{dz}\Big)f(a_lz)=-w_1\Big(z\frac{d}{dz}\Big)f(a_1z)-\ldots-w_{l-1}\Big(z\frac{d}{dz}\Big)f(a_{l-1}z)+\frac{1}{1-z}.$$
After applying the equality above inductively, we receive for any $N\in\NN$
\begin{multline*}
w_l\Big(z\frac{d}{dz}\Big)^Nf(a_l^Nz)=w_l\Big(z\frac{d}{dz}\Big)^{N-1}\frac{1}{1-a_l^{N-1}z}-w_1\Big(z\frac{d}{dz}\Big)w_l\Big(z\frac{d}{dz}\Big)^{N-1}f(a_l^{N-1}a_1z)-\\-\ldots-w_{l-1}\Big(z\frac{d}{dz}\Big)w_l\Big(z\frac{d}{dz}\Big)f(a_l^{N-1}a_{l-1}z)=\ldots=\\ =\sum_{K=1}^N(-1)^{K-1}w_l\Big(z\frac{d}{dz}\Big)^{N-K}\bigg(\sum_{i_1+\ldots+i_{l-1}=K-1}\frac{(K-1)!}{i_1!\ldots i_{l-1}!}w_1\Big(z\frac{d}{dz}\Big)^{i_1}\ldots\\\ldots w_{l-1}\Big(z\frac{d}{dz}\Big)^{i_{l-1}}\frac{1}{1-a_1^{i_1}\ldots a_{l-1}^{i_{l-1}}a_l^{N-K}z }\bigg)+ \\ +(-1)^N\sum_{i_1+\ldots+i_{l-1}=N}\frac{N!}{i_1!\ldots i_{l-1}!}w_1\Big(z\frac{d}{dz}\Big)^{i_1}\ldots w_{l-1}\Big(z\frac{d}{dz}\Big)^{i_{l-1}}f\Big(a_1^{i_1}\ldots a_{l-1}^{i_{l-1}}z\Big).
\end{multline*}

Hence,
\begin{multline}
\label{thm (9)}
f(a_l^Nz)=\sum_{K=1}^N(-1)^{K-1}\Bigg(\sum_{i_1+\ldots+i_{l-1}=K-1}\frac{(K-1)!}{i_1!\ldots i_{l-1}!}\bigg[\frac{w_1\Big(z\frac{d}{dz}\Big)}{w_l\Big(z\frac{d}{dz}\Big)}\bigg]^{i_1}\ldots \\\ldots\bigg[\frac{w_{l-1}\Big(z\frac{d}{dz}\Big)}{w_l\Big(z\frac{d}{dz}\Big)}\bigg]^{i_{l-1}}\frac{1}{1-a_1^{i_1}\ldots a_{l-1}^{i_{l-1}}a_l^{N-K}z }\Bigg)+\\+(-1)^N\sum_{i_1+\ldots+i_{l-1}=N}\bigg[\frac{w_1\Big(z\frac{d}{dz}\Big)}{w_l\Big(z\frac{d}{dz}\Big)}\bigg]^{i_1}\ldots \bigg[\frac{w_{l-1}\Big(z\frac{d}{dz}\Big)}{w_l\Big(z\frac{d}{dz}\Big)}\bigg]^{i_{l-1}}f\Big(a_1^{i_1}\ldots a_{l-1}^{i_{l-1}}z\Big),
\end{multline}
where the linear operator $\displaystyle\frac{\tilde w_1\Big(z\frac{d}{dz}\Big)}{\tilde w_2\Big(z\frac{d}{dz}\Big)}:\CC[[z]]\longrightarrow\CC[[z]]$ is defined as $$\frac{\tilde w_1\Big(z\frac{d}{dz}\Big)}{\tilde w_2\Big(z\frac{d}{dz}\Big)}\bigg(\sum_{n=0}^\infty a_nz^n\bigg):=\sum_{n=0}^\infty \frac{\tilde w_1(n)}{\tilde w_2(n)}a_nz^n.$$
Formula~(\ref{thm (9)}) gives analytic continuation of the function $f\in\Oo(D_{a_l})$ to the set $\CC\setminus[a_l,\infty)$. To estimate $f(z)$, first observe that for any $z\in L_{\theta}$, $\theta\neq 0\mod 2\pi$, and for $N=\left\lceil\log_{\left(\frac{a_l}{a_{l-1}}\right)}\left(\frac{2|z|}{a_l}\right)\right\rceil$ and $i_1+\ldots+i_{l-1}=N$ we have
$$ \Bigg|\frac{a_1^{i_1}\ldots a_{l-1}^{i_{l-1}}z}{a_l^N}\Bigg|\leq \Big(\frac{a_{l-1}}{a_l}\Big)^N|z|\leq \frac{a_l}{2}. $$
Hence $$ \Bigg|f\bigg(\frac{a_1^{i_1}\ldots a_{l-1}^{i_{l-1}}z}{a_l^N}\bigg)\Bigg|\leq\mathrm{sup}_{|\zeta|=\frac{a_l}{2}}\Big|f(\zeta)\Big|<\infty.$$
Moreover, by Cauchy's estimates there exist constants $A,B>0$ and $p\in\NN$ independent of $N$ such that
\begin{multline*}
\left|(-1)^N\sum_{i_1+\ldots+i_{l-1}=N}\frac{N!}{i_1!\ldots i_{l-1}!}\left[\frac{w_1\Big(z\frac{d}{dz}\Big)}{w_l\Big(z\frac{d}{dz}\Big)}\right]^{i_1}\ldots \left[\frac{w_{l-1}\Big(z\frac{d}{dz}\Big)}{w_l\Big(z\frac{d}{dz}\Big)}\right]^{i_{l-1}}f\bigg(\frac{a_1^{i_1}\ldots a_{l-1}^{i_{l-1}}z}{a_l^N}\bigg)\right|\\ \leq AB^N(N!)^p.
\end{multline*}
Since $$N\leq\log_{\big(\frac{a_l}{a_{l-1}}\big)}\Big(\frac{2|z|}{a_l}\Big)+1=\frac{\ln|z|+\ln\Big(\frac{2}{a_l}\Big)}{\ln\Big(\frac{a_l}{a_{l-1}}\Big)}+1,$$
we observe that for sufficiently large $|z|$ there exists a constant $d>0$ such that $N\leq d\ln|z|$. It means that 
\begin{equation}\label{thm (10)}
AB^N(N!)^p\leq AB^{d\ln|z|}\Big(d\ln|z|\Big)^{pd\ln|z|}\leq \tilde Ae^{\tilde B \ln|z|\Big(1+\ln\big(\ln|z|\big)\Big)}
\end{equation}
for some constants $\tilde A,\tilde B<\infty$.
Analogously, from Lemma \ref{lemat 3} follows the existence of $C<\infty$, $p\in\NN$ and $n_0\in\NN$ such that
\begin{multline*}
\Bigg|\sum_{K=1}^N(-1)^{K-1}\bigg(\sum_{i_1+\ldots+i_{l-1}=K-1}\frac{(K-1)!}{i_1!\ldots i_{l-1}!}\bigg[\frac{w_1\Big(z\frac{d}{dz}\Big)}{w_l\Big(z\frac{d}{dz}\Big)}\bigg]^{i_1}\ldots\\\ldots\bigg[\frac{w_{l-1}\Big(z\frac{d}{dz}\Big)}{w_l\Big(z\frac{d}{dz}\Big)}\bigg]^{i_{l-1}}\frac{1}{1-a_1^{i_1}\ldots a_{l-1}^{i_{l-1}}a_l^{-K}z }\bigg)\Bigg|\leq\sum_{K=1}^N\frac{C^K\big((Kp)!\big)^2}{C(\theta)^{pK+1}}\Big(|z|^{n_0}+1\Big)\leq\\\leq\tilde C^N\big((Np)!\big)^2\Big(|z|^{n_0}+1\Big).
\end{multline*}
For sufficiently large $|z|$, using inequality $N\leq d\ln|z|$ we may estimate it by 
\begin{equation}\label{thm (11)}
\tilde C^{d\ln|z|}\Big(dp\ln|z|\Big)^{2dp\ln|z|}|z|^{n_0}\leq A^{'}e^{B^{'}\ln|z|\Big(1+\ln\big(\ln|z|\big)\Big)}
\end{equation}
for some constants $A^{'},B^{'}<\infty$.
Using~(\ref{thm (10)}) and~(\ref{thm (11)}) we observe the fact that $f(z)\in\Oo^{>0}(\CC\setminus\RR_+)$. Hence,
by Theorem~\ref{th:1} we conclude that the sequence $(m(n))_{n\geq 0}$ preserves summability. 
\end{proof}

Since the set of sequences preserving summability forms a group with multiplication (see Remark \ref{re:groups}),
we may replace polynomials given in Theorem \ref{th:10} by rational functions. In this way we receive the most general class of sequences satisfying summability among considered in this section.
\begin{Cor}
Let 
\begin{equation*}
m(n)=\frac{w_1(n)}{v_1(n)}a_1^n+\ldots+\frac{w_l(n)}{v_l(n)}a_l^n\quad\textrm{for every}\quad n\in\NN_0
\end{equation*}
be a sequence of positive numbers, where $0<a_1<\ldots<a_l$ and $w_1,\ldots,w_l,v_1,\ldots,v_l$ are given polynomials with complex coefficients,  $w_l(n)\neq0$ for every $n\in\NN_0$, and $v_i(n)\neq 0$ for every $n\in\NN_0$ and $i=1,\ldots,l$. Then the sequence $(m(n))_{n\in\NN_0}$ preserves summability.
\end{Cor}

\section{Applications and final remarks}\label{sec:5}
\subsection{Applications to moment differential equations}
In this section we apply sequences preserving summability to the study of summable solutions of some moment differential equations.

The notion of $m$-moment differentiation was introduced by Balser and Yoshino~\cite{B-Y} in the case when $m$ is a moment sequence associated with a pair of kernel functions (see \cite[Section 5.5]{B2}) and then extended in~\cite{I-M} to any sequence $m=(m(n))_{n\geq 0}$ of positive numbers with $m(0)=1$. 
The solutions of moment differential equations where studied in such papers as \cite{L-M-S2, L-M-S4, L-P, L-Sa-Se, Mic7, Sa}.
\begin{Def}
 For a given sequence $m=(m(n))_{n\geq 0}$ of positive numbers with $m(0)=1$, an operator $\partial_{m,t}\colon\EE[[t]] \to \EE[[t]]$ defined by
 \begin{equation*}
  \partial_{m,t}\big(\sum_{n=0}^{\infty}u_nt^n\big):=\sum_{n=0}^{\infty}\frac{m(n+1)}{m(n)}u_{n+1}t^n
 \end{equation*}
is called an \emph{$m$-moment differential operator}.
\end{Def}

\begin{Rem}
 Observe that in the most important case $\Gamma_1:=(\Gamma(1+n))_{n\geq 0}=(n!)_{n\geq 0}$, the operator $\partial_{\Gamma_1,t}$ is the $\Gamma_1$-moment differential operator, which coincides with the usual derivative $\partial_t$.
\end{Rem}

By the direct calculation we receive the following commutation formula between $m_1$-Borel operator and $m_2$-moment differentiation.
\begin{Prop}
\label{prop:commutate}
 Let $m_1=(m_1(n))_{n\geq 0}$ and $m_2=(m_2(n))_{n\geq 0}$ be sequences of positive numbers. Then the operators $\Bo_{m_1},\partial_{m_2,t}\colon\EE[[t]]\to\EE[[t]]$ commute in a such way that
 \begin{equation*}
  \Bo_{m_1}\partial_{m_2,t}=\partial_{m_1m_2,t}\Bo_{m_1}.
 \end{equation*}
\end{Prop}

We assume that $m=(m(n))_{n\geq 0}$ is a sequence of positive numbers with $m(0)=1$  and $P(\lambda,\zeta)$ is a general polynomial of two variables of order $p$ with respect to $\lambda$ and $\varphi_j(z)\in\Oo(D)$ for $j=0,\dots,p-1$.

The focus of our study is the relationship between the solution $\hat{u}(t,z)\in\Oo(D)[[t]]$ of the Cauchy problem
\begin{equation}
\label{eq:CP_u}
    \left\{
    \begin{aligned}
     P(\partial_{m\Gamma_1,t},\partial_z)u&=0\\
     \partial^j_{m\Gamma_1,t}u(0,z)&=\varphi_j(z),\quad j=0,\dots,p-1,
    \end{aligned}
    \right.
 \end{equation}
and the solution $\hat{v}(t,z)\in\Oo(D)[[t]]$ of a similar initial value problem of the form
\begin{equation}
 \label{eq:CP_tilde_u}
    \left\{
    \begin{aligned}
     P(\partial_{t},\partial_z)v&=0\\
     \partial^j_{t}v(0,z)&=\varphi_j(z),\quad j=0,\dots,p-1.
    \end{aligned}
    \right.
 \end{equation}

 Analogously to \cite[Proposition 7]{I-M}, where the similarity between the operators $\partial_{m,t}$ and $\partial_{\mathbf{1},t}$, both of order zero, was investigated, we study the connection between the operators of order one, $\partial_{m\Gamma_1,t}$ and $\partial_{t}$. 
 We have (see also \cite[Proposition 7]{I-M})
 \begin{Prop}\label{prop:9}
  Let $m=(m(n))_{n\geq 0}$ be a sequence of positive numbers with $m(0)=1$. Then $\hat{u}(t,z)=\sum_{n=0}^{\infty}\frac{u_n(z)}{m(n)}t^n$ is a formal power series solution of (\ref{eq:CP_u}) if and only if $\hat{v}(t,z)=\sum_{n=0}^{\infty}u_n(z)t^n$ is a formal power series solution of (\ref{eq:CP_tilde_u}).
 \end{Prop}
\begin{proof}
 ($\Rightarrow$) Let $\hat{u}(t,z)=\sum_{n=0}^{\infty}\frac{u_n(z)}{m(n)}t^n$ be a formal solution of (\ref{eq:CP_u}). Using the commutation formula (Proposition \ref{prop:commutate}) $\Bo_{m^{-1}}\partial_{m\Gamma_1,t} = \partial_{t}\Bo_{m^{-1}}$  and applying the Borel transform $\Bo_{m^{-1}}$ to the Cauchy problem (\ref{eq:CP_u}) we conclude that $\hat{v}(t,z)=\Bo_{m^{-1}}\hat{u}(t,z)$ is a formal solution of (\ref{eq:CP_tilde_u}).
\medskip\par
($\Leftarrow$) The proof is analogous. It is sufficient to apply the Borel transform $\Bo_{m,t}$ to the Cauchy problem (\ref{eq:CP_tilde_u}) and to observe that 
$\Bo_{m}\partial_{t} = \partial_{m\Gamma_1,t}\Bo_{m}$ and
$\hat{u}(t,z)=\Bo_{m}\hat{v}(t,z)$.
\end{proof}

Directly from Proposition~\ref{prop:9} and Definition~\ref{df:preserve} we receive
\begin{Prop}
Let $m=(m(n))_{n\geq 0}$ be a sequence preserving summability, $k>0$ and $d\in\RR$. Then a formal power series solution $\hat{u}$ of (\ref{eq:CP_u}) is $k$-summable in a direction $d$ if and only if a formal power series solution $\hat{v}$ of (\ref{eq:CP_tilde_u}) is $k$-summable in the same direction $d$.  
\end{Prop}

It means that the moment differential operator $\partial_{m\Gamma_1}$ has the same behavior in the context of summability as the usual differential operator $\partial_t$. In particular, we receive the following version of the result \cite[Theorem 3.1]{Lu-Mi-Sc} by Lutz, Miyake and Sch\"afke that can be applied to summable solutions of the heat equation:
\begin{Cor}
Let us assume that $m=(m(n))_{n\geq 0}$ is a sequence preserving summability, $d\in\RR$ and $\hat{u}\in\Oo(D)[[t]]$ is a formal solution of the Cauchy problem
\begin{equation*}
    \left\{
    \begin{aligned}
     \partial_{m\Gamma_1,t}u - \partial_z^2u&=0\\
     u(0,z)&=\varphi(z)\in\Oo(D).
    \end{aligned}
    \right.
 \end{equation*}
Then $\hat{u}$ is $1$-summable in a direction $d$ if and only if $\varphi\in\Oo^2(\hat{S}_{d/2}\cup \hat{S}_{d/2+\pi})$.    
\end{Cor}

\begin{Ex}
The following operators have the same behavior in the context of summability as the usual differential operator $\partial_t$.
\begin{enumerate}
    \item[1.] $4\partial_t-2\partial_{\mathbf{1},t}=\partial_{m\Gamma_1,t}$ for $m=\left(\frac{(2n)!}{(n!)^2}\right)_{n\geq 0}$ being a sequence of moments of order zero given in Example \ref{ex:4}.
    \item[2.] $\partial_t D_{q,t}=\partial_{m\Gamma_1,t}$ for
    $m=([n]_q!)$, where $q\in[0,1)$. Here $D_{q,t}$ denotes the $q$-difference operator defined as $D_{q,t}u(t)=\frac{u(qt)-u(t)}{qt-t}$ (see Example \ref{ex:5}).
    \item[3.] $a\partial_t(1+\frac{1}{t}\partial_t^{-1})^p=\partial_{m\Gamma_1,t}$ for $m=\left(a^nw(n)\right)_{n\geq 0}=\left(a^n(n+2)^p)\right)_{n\geq 0}$, where $a>0$, $p\in\NN$ and $\partial_t^{-1}$ denotes the antiderivative defined as $\partial_t^{-1}u(t)=\int_0^tu(\tau)\,d\tau$. Observe that the sequence $m$ preserves summability by Theorem \ref{th:10}.  
\end{enumerate}
\end{Ex}

\subsection{Final remarks}
Directly from Definition~\ref{df:preserve} it follows that the product of sequences preserving summability also preserves summability (see also Remark \ref{re:groups}). It is worth asking whether the sum of sequences preserving summability also preserves summability. Condition $m(0)=1$ forces the need to modify the definition of a sum of sequences a little bit. Namely, for sequences $m_1=(m_1(n))_{n\geq 0}$ and $m_2=(m_2(n))_{n\geq 0}$ we define their sum as
\begin{equation}\label{eq:sum}
m_1+m_2=(m(n))_{n\geq0},\quad\textrm{where}\quad
m(n)=\left\{
    \begin{array}{lll}
     1&\textrm{for}&n=0\\
     m_1(n)+m_2(n)&\textrm{for}&n\geq 1
    \end{array}
    \right..
\end{equation}
By Theorem \ref{th:1} and Remark \ref{re:4}, to prove that $m$ is a sequence preserving summability it is sufficient to show that
\begin{equation}\label{eq:m-}
\Bo_{m^{-1}}\left(\sum_{n=0}^{\infty} t^n\right)\in\Oo^{>0}(\CC\setminus\RR_+)
\end{equation}
and 
\begin{equation}\label{eq:m+}
\Bo_{m}\left(\sum_{n=0}^{\infty} t^n\right)\in\Oo^{>0}(\CC\setminus\RR_+)
\end{equation}
Since 
$$
\Bo_{m^{-1}}\left(\sum_{n=0}^{\infty} t^n\right)=\Bo_{m_1^{-1}}\left(\sum_{n=0}^{\infty} t^n\right) + \Bo_{m_2^{-1}}\left(\sum_{n=0}^{\infty} t^n\right) - 1
$$
and $m_1$, $m_2$ preserve summability, by Theorem \ref{th:1} and Remark \ref{re:4}
we conclude that $\Bo_{m_1^{-1}}(\sum t^n)\in\Oo^{>0}(\CC\setminus\RR_+)$ and $\Bo_{m_2^{-1}}(\sum t^n)\in\Oo^{>0}(\CC\setminus\RR_+)$, and
also \eqref{eq:m-} holds.

Unfortunately, it is not possible to express
in a similar way $\Bo_{m}(\sum t^n)$ in terms of 
$\Bo_{m_1}(\sum t^n)$, $\Bo_{m_2}(\sum t^n)$, $\Bo_{m_1^{-1}}(\sum t^n)$ and $\Bo_{m_2^{-1}}(\sum t^n)$, so we do not know if \eqref{eq:m+} also holds. Moreover, we suppose that the set of sequences preserving summability is not closed under addition, but we have no proof of this conjecture.


\begin{thebibliography}{30}


\bibitem{B2} W. Balser, \textit{Formal power series and linear systems of meromorphic ordinary differential equations}, Springer-Verlag, New York, 2000.

%
\bibitem{B-Y} W. Balser, M. Yoshino, \textit{Gevrey order of formal power series solutions of inhomogeneous partial differential equations with constant coefficients}, Funkcial. Ekvac. 53 (2010), 411--434.

\bibitem{E-M-O-T} A. Erdélyi, W. Magnus, F. Oberhettinger, F.G. Tricomi, \emph{Higher Transcendental Functions, vol. I},
McGraw–Hill, New York, 1953.

\bibitem{G-K-P} R. Graham, D. Knuth, O. Patashnik, \emph{Concrete Mathematics. A foundation for computer science}, 2nd ed.,  Addison-Wesley Publishing Company, Reading, 1994.

\bibitem{G-S} J. Guillera, J. Sondow, \emph{Double integrals and inﬁnite products for some classical constants via analytic continuations of Lerch’s transcendent}, Ramanujan J. 16 (2008) 247–270.
%
\bibitem{I-M} K. Ichinobe, S. Michalik, \textit{On the summability and convergence of formal solutions of linear $q$-difference-differential equations with constant coefficients}, Math. Ann. 389 (2024), 1099--1130.

\bibitem{L-M} A. Lastra, S. Michalik, \textit{On sequences preserving $q$-Gevrey asymptotic expansions}, Anal. Math. Phys. 14, 17 (2024).

\bibitem{L-M-S2} A. Lastra, S. Michalik, M. Suwi\'nska, \emph{Summability of formal solutions for some generalized moment partial differential equations}, Results Math. 76:22 (2021).

\bibitem{L-M-S4} A. Lastra, S. Michalik, M. Suwi\'nska, \emph{Multisummability of formal solutions for a family of generalized singularly perturbed moment differential equations}, Results Math 78:49 (2023).

\bibitem{L-P} A. Lastra, C. Prisuelos-Arribas,
\emph{Solutions of linear systems of moment differential equations via generalized matrix exponentials}, J.~Differ. Equations 372 (2023), 591--611.

\bibitem{L-Sa-Se} A. Lastra, J.~Sanz, J.~R. Sendra,
\emph{On the summability of a class of formal power series}, Math. Inequal. Appl. 25 (2022), 1101--1121.


%
%
\bibitem{LR} M. Loday-Richaud, \textit{Divergent series, summability and resurgence II. Simple and multiple summability}, Lecture Notes in Mathematics, vol. 2154, Springer 2016.

\bibitem{Lu-Mi-Sc} D. A. Lutz, M. Miyake, R. Scha\"afke, \emph{On the Borel summability of divergent solutions of the heat equation}, Nagoya Math. J. 154 (1999), 1--29.

%
%
%
\bibitem{Mic7} S. Michalik, \textit{Analytic solutions of moment partial differential equations with constant coefficients}, Funkcial. Ekvac. 56 (2013), 19--50.
%
\bibitem{Mic8} S. Michalik, \textit{Summability of formal solutions of linear partial differential equations with divergent initial data}, J. Math. Anal. Appl. 406 (2013), 243--260.
%
%
%
\bibitem{Sa} J.~Sanz, \emph{Asymptotic analysis and summability of formal power series}, Analytic, algebraic and geometric aspects of differential equations, Trends in Mathematics (2017), 199--262.
%
%
%


\end{thebibliography}
\end{document}